\documentclass[a4paper,leqno,draft]{article}
\usepackage{amssymb,amsmath,amsfonts,amsthm,
mathrsfs}

\renewcommand{\theequation}{\thesection.\arabic{equation}}
\newtheorem{theorem}{Theorem}[section]
\newtheorem{corollary}[theorem]{Corollary}
\newtheorem{lemma}[theorem]{Lemma}
\newtheorem{proposition}[theorem]{Proposition}
\newtheorem{definition}[theorem]{Definition}
\theoremstyle{definition}
\newtheorem{remark}[theorem]{Remark}

\newcommand\dslash{d\llap {\raisebox{.9ex}{$\scriptstyle-\!$}}}
\newcommand{\wt}[1]{\widetilde{#1}}

\newcommand{\Cinf}{\ensuremath{\mathcal{C}^\infty}}
\newcommand{\Cinfc}{\ensuremath{\mathcal{C}^\infty_{\text{c}}}}
\newcommand{\D}{\ensuremath{{\cal D}}}
\renewcommand{\S}{\mathscr{S}}
\newcommand{\E}{\ensuremath{{\cal E}}}

\newcommand{\mb}[1]{\ensuremath{\mathbb{#1}}}
\newcommand{\N}{\mb{N}}

\newcommand{\R}{\mb{R}}
\newcommand{\C}{\mb{C}}

\newcommand{\cl}[1]{\ensuremath{[#1]}}

\newcommand{\G}{\ensuremath{{\cal G}}}

\newcommand{\Gc}{\ensuremath{{\cal G}_\mathrm{c}}}
\newcommand{\Gcinf}{\ensuremath{{\cal G}^\infty_\mathrm{c}}}

\newcommand{\EM}{\ensuremath{{\cal E}_{M}}}

\newcommand{\EMinf}{\ensuremath{{\cal E}^\infty_{M}}}
\newcommand{\ESinf}{\mathcal{E}_{\S}^{\infty}}

\newcommand{\NN}{\ensuremath{{\cal N}}}

\newcommand{\g}{\mathrm{g}} 
\newcommand{\Ginf}{\ensuremath{\G^\infty}}

\newcommand{\lara}[1]{\langle #1 \rangle}
\newcommand{\WF}{\mathrm{WF}}
\newcommand{\WFg}{\WF_{\mathrm{g}}}
\newcommand{\singsupp}{\mathrm{sing supp}}
\newcommand{\supp}{\mathrm{supp}}
\newcommand{\Char}{\ensuremath{\text{Char}}}
\newcommand{\zs}{\setminus 0}
\newcommand{\CO}[1]{\ensuremath{T^*(#1) \zs}}
\newcommand{\ssc}{\mathrm{sc}}

\newcommand{\Ellsc}{\mathrm{Ell}_\ssc}

\newfont{\bigmath}{cmr12 at 13pt}
\newcommand{\PPsi}{\bigmath{\symbol{9}}}
\newcommand{\Oprop}[1]{{\ }_{\mathrm{pr}}^{\hphantom{m}}\text{\PPsi}_{\ssc}^{ #1}}
\newcommand{\Opropr}[1]{{\ }_{\mathrm{pr}}^{\hphantom{m}}\text{\PPsi}_{\rm{rg}}^{ #1}}

\newfont{\grecomath}{cmmi12 at 15pt}
\newcommand{\Bigmu}{\text{\grecomath{\symbol{22}}}}

\renewcommand{\d}{\ensuremath{\partial}}

\newfont{\bl}{msbm10 scaled \magstep2}

\newcommand{\beq}{\begin{equation}}
\newcommand{\eeq}{\end{equation}}

\newcommand{\emb}{\hookrightarrow}

\newcommand{\FT}[1]{\widehat{#1}}

\newcommand{\F}{\ensuremath{{\cal F}}}

\newcommand{\notmid}{\mid\kern-0.5em\not\kern0.5em}

\newcommand{\norm}[2]{{\| #1 \|}_{#2}}

\newcommand{\al}{\alpha}

\newcommand{\ga}{\gamma}
\newcommand{\Ga}{\Gamma}

\newcommand{\eps}{\varepsilon}

\newcommand{\la}{\lambda}

\newcommand{\Om}{\Omega}

\newcommand{\compl}[1]{{#1}^{\mathrm{c}}}

\newcommand{\Sscu}{\underline{\mathcal{S}}_{\,\ssc}}
\newcommand{\Nu}{\underline{\mathcal{N}}}
\newcommand{\Syscu}{{\wt{\underline{\mathcal{S}}}}_{\,\ssc}}
\newcommand{\Nuinf}{\Nu^{-\infty}}
\newcommand{\Sbeu}{\underline{\mathcal{S}}_{\mathfrak{B}}}
\newcommand{\Sybeu}{{\wt{\underline{\mathcal{S}}}}_{\mathfrak{B}}}
\newcommand{\Syru}{{\wt{\underline{\mathcal{S}}}}_{\mathrm{rg}}}
\newcommand{\Sru}{{{\underline{\mathcal{S}}}}_{\mathrm{rg}}}

\newcommand{\Syr}{{\wt{{\mathcal{S}}}}_{\mathrm{rg}}}

\newcommand{\Wsc}{\mathrm{W}_{\ssc}}
\newcommand{\Wcl}{\mathrm{W}_{\rm{cl}}}
\newcommand{\mufty}{\mu_\g}

\renewcommand{\dslash}{d\hspace{-0.4em}{ }^-\hspace{-0.2em}}

\begin{document}

\title{{\bf Microlocal analysis of generalized functions:
pseudodifferential techniques and propagation of singularities}}

\author{Claudia Garetto \\
Dipartimento di Matematica\\ Universit\`a di Torino, Italia\\
\texttt{garettoc@dm.unito.it}\\
\ \\
G\"{u}nther H\"{o}rmann\footnote{Supported by FWF grant P14576-MAT,
current address: Institut f\"ur Technische Mathematik, Geometrie und Bauniformatik,
Universit\"at Innsbruck}\\
Institut f\"ur Mathematik\\
Universit\"at Wien, Austria\\
\texttt{guenther.hoermann@univie.ac.at}
}
\date{February 2004}
\maketitle

\begin{abstract} We characterize microlocal regularity, in the $\Ginf$-sense,
of Colombeau generalized functions by an appropriate extension of the classical
notion of micro-ellipticity to pseudodifferential operators with slow scale
generalized symbols. Thus we obtain an alternative, yet equivalent, way to
determine generalized wave front sets, which is analogous to the original
definition of the wave front set of distributions via intersections over
characteristic sets. The new methods are then applied to regularity theory of
generalized solutions of (pseudo-)differential equations, where we extend the
general noncharacteristic regularity result for distributional solutions and
consider propagation of $\Ginf$-singularities for homogeneous first-order
hyperbolic equations. 
\end{abstract}

\setcounter{section}{-1}
\section{Introduction}

Microlocal analysis in Colombeau algebras of generalized functions, as it has
been initiated (in published form) in \cite{DPS:98,NPS:98}, is a compatible
extension of its distribution theoretic analogue to the realm of an
unrestricted differential-algebraic context. The main emphasis in recent
research on the subject is on microlocal properties of basic nonlinear
operations as well as on regularity theory for generalized solutions to partial
(and pseudo-) differential equations (cf.\ \cite{NP:98,HK:01,HdH:01,HOP:03,GH:03,Garetto:04,GGO:03}).

For Schwartz distributions the so-called elementary characterization of
microlocal regularity is a corollary to its original definition via
characteristic sets under pseudodifferential actions (cf.\cite{Hoermander:71}),
whereas the intuitively appealing function-like aspect of Colombeau generalized
functions seems to have fostered a ``generalized elementary'' approach as being
a natural definition there. This may have two reasons: first, the new
microlocal regularity notion is based on $\Ginf$-regularity, which coincides
with $\Cinf$-regularity in case of embedded distributions, and was introduced
in 
\cite{O:92} in direct analytical terms in form of asymptotic estimates of the
derivatives; second, as soon as one puts oneself into the much wider setting of
Colombeau spaces -- with the possibility to allow for highly singular symbols
of (pseudo-)differential operators -- the question of good choices for
appropriate generalized notions of the characteristic set or
(micro-)ellipticity turns into a considerable and crucial part of the research
issue (cf.\ \cite{HdH:01,HO:03,HOP:03,GGO:03}).

In the present paper, we succeed to prove characterizations of the generalized
wave front set of a Colombeau generalized functions in terms of intersections
over certain non-ellipticity domains corresponding to pseudodifferential
operators yielding $\Ginf$-regular images. Thus we obtain direct analogues of
H\"ormander's definition of the distributional wave front set given in
\cite{Hoermander:71}. Moreover, as first test applications of the new results
we discuss a generalization of the noncharacteristic regularity theorem for
pseudodifferential equations and propagation of $\Ginf$-singularities (or
rather, generalized wave front sets) for generalized solutions of first-order
hyperbolic differential equations with $\Cinf$-coefficients.

In the remainder of the introductory section we fix some notation and review
basic notions from Colombeau theory. Section 1 provides the technical
background on the generalized symbol classes used later on and introduces an
appropriate micro-ellipticity notion. The theoretical core of the paper is
Section 2 where the main results on micro-locality, microsupport, and the wave
front set characterizations are proven. Section 3 discusses applications to
regularity theory of generalized solutions of (pseudo-)differential equations.
Since by now several variants of pseudodifferential operator approaches in
Colombeau algebras and generalized symbol calculi occur in the literature
(\cite{NPS:98,Garetto:04,GGO:03}), and we employ yet a slightly different one
here, we decided to sketch the basics of such a rather general scheme of
calculus in the Appendix, the skeleton of which is structurally close to the
comprehensive treatment in \cite{GGO:03}.

We point out that \cite{NPS:98} already includes a result on micro-locality
(similar to our Theorem \ref{theorem_micro_support}) of actions of generalized
pseudodifferential operators, whose definition is based solely on regularizing
nets of symbols, but not Colombeau classes, and uses Fourier integral
representations with additional asymptotic cut-offs. The definitions of the
operator actions can thus be compared in a weak sense only.

\subsection{Notation and basic notions from Colombeau theory}

We use \cite{Colombeau:84,Colombeau:85,O:92,GKOS:01} as standard references for
the foundations and various applications of Colombeau theory. We will work with
the so-called special Colombeau algebras, denoted by $\G^s$ in \cite{GKOS:01},
although here we will consistently drop the superscript $s$ to avoid notational
overload.

We briefly recall the basic construction. Throughout the paper $\Om$ will
denote an open subset of $\R^n$. \emph{Colombeau generalized functions} on
$\Om$ are defined as equivalence classes $u = \cl{(u_\eps)_\eps}$ of nets of
smooth functions $u_\eps\in\Cinf(\Om)$ (\emph{regularizations}) subjected to
asymptotic norm conditions with respect to $\eps\in (0,1]$ for their
derivatives on compact sets. More precisely, we have
\begin{itemize}
 \item moderate nets $\EM(\Om)$: $(u_\eps)_\eps\in\Cinf(\Om)^{(0,1]}$ such that
 for all $K \Subset \Om$ and $\al\in\N^n$ there exists $p \in \R$ such that
 \beq\label{basic_estimate}
    \norm{\d^\al u_\eps}{L^\infty(K)} = O(\eps^{-p}) \qquad (\eps \to 0);
 \eeq
  \item negligible nets $\NN(\Om)$: $(u_\eps)_\eps\in \EM(\Om)$ such that for all
  $K \Subset \Om$ and for all $q\in\R$ an
 estimate $\norm{u_\eps}{L^\infty(K)} = O(\eps^{q})$ ($\eps \to 0$) holds;
 \item $\EM(\Om)$ is a differential algebra with operations defined at fixed $\eps$,
 $\NN(\Om)$ is an ideal, and $\G(\Om) := \EM(\Om) / \NN(\Om)$ is the (special)
 \emph{Colombeau algebra};
 \item there are embeddings, $\Cinf(\Om)\emb \G(\Om)$ as subalgebra and
 $\D'(\Om) \emb \G(\Om)$ as linear space, commuting with partial derivatives;
 \item  $\Om \to \G(\Om)$ is a fine sheaf and $\Gc(\Om)$ denotes
 the subalgebra of elements with compact support; by a
 cut-off in a neighborhood of the support one can always obtain
 representing nets with supports contained in a joint compact set.
\end{itemize}

The subalgebra $\Ginf(\Om)$ of \emph{regular Colombeau}, or
\emph{$\Ginf$-regular}, generalized functions consists of those elements in
$\G(\Om)$ possessing representatives such that estimate (\ref{basic_estimate})
holds for a certain $m$ uniformly over all $\al\in\N^n$. We will occasionally
use the notation $\EMinf(\Om)$ for the set of such nets of regularizations. In
a similar sense, we will denote by $\ESinf(\R^n)$ the set of regularizations
$(v_\eps)_\eps\in\S(\R^n)^{(0,1]}$ with a uniform asymptotic power of $\eps$-growth for
all $\S$-seminorms of $v_\eps$.

A Colombeau generalized function $u = [(u_\eps)_\eps]\in\G(\Om)$ is said to be
\emph{generalized microlocally regular}, or \emph{$\Ginf$-microlocally
regular}, at $(x_0,\xi_0)\in \CO{\Om} = \Om\times (\R^n\setminus\{0\})$
(cotangent bundle with the zero section removed) if there is 
$\phi\in\Cinfc(\Om)$
with $\phi(x_0) = 1$ and a conic neighborhood $\Ga\subseteq\R^n\setminus\{0\}$
of $\xi_0$ such that $\F(\phi u)$ is \emph{(Colombeau) rapidly
decreasing} in $\Ga$ (cf.\ \cite{Hoermann:99}), i.e., there exists $N$ such that for all $l$ we have
 \beq\label{rap_dec_Ga}
 \sup\limits_{\xi\in\Ga} \lara{\xi}^l |(\phi u_\eps)\FT{\ }(\xi)| = O(\eps^{-N})
    \qquad (\eps \to 0),
\eeq where we have used the standard notation $\lara{\xi} = (1 +
|\xi|^2)^{1/2}$. Note that, instead of specifying the test function $\phi$ as
above, one may equivalently require the existence of an open neighborhood $U$
of $x_0$ such that for all $\phi\in\Cinfc(U)$ the estimate (\ref{rap_dec_Ga})
holds.  

Finally, we will use the term \emph{proper cut-off function} for any
$\chi\in\Cinf(\Om\times\Om)$ such that $\supp(\chi)$ is a proper subset of
$\Om\times\Om$ (i.e., both projections are proper maps) and $\chi =1$ in a
neighborhood of the diagonal $\{(x,x) : x\in\Om\}\subset\Om\times\Om$.

\section{Slow scale micro-ellipticity}

The pseudodifferential operator techniques which we employ are based on a generalization of the classical symbols spaces $S^m(\Om\times\R^n)$ (cf.\ \cite{Hoermander:71}). These spaces are Fr\'{e}chet spaces endowed with the seminorms 
\beq\label{seminorm}
 |a|^{(m)}_{K,\alpha,\beta} = \sup_{x\in K,\xi\in\R^n}\langle\xi\rangle^{-m+|\alpha|}|\partial^\alpha_\xi\partial^\beta_x a(x,\xi)|,
\eeq
where $K$ ranges over the compact subsets of $\Om$.  Several types of Colombeau generalized symbols are studied in \cite{NPS:98, Garetto:04, GGO:03, GH:03} providing a pseudodifferential calculus and regularity theory on the level of operators. In the current paper, we focus on the microlocal aspects. To this end we introduce a particularly flexible class of symbols with good stability properties with respect to lower order perturbations. In the spirit of the earlier Colombeau approaches, our symbols are defined via families $(a_\eps)_{\eps}\in S^m(\Om\times\R^n)^{(0,1]} =: {\mathcal{S}}^m[\Om\times\R^n]$ of regularizations, subjected to asymptotic estimates of the above seminorms in terms of $\eps$. The particular new feature is a slow scale growth, which proved to be essential in regularity theory (cf. \cite{HO:03, HOP:03, GGO:03}). This property is measured by the elements of the following set of \emph{strongly positive slow scale nets}
\beq\label{slow_scale}\begin{split}
  \Pi_\ssc := \{ (\omega_\eps)_\eps\in\R^{(0,1]}\, :\quad &\exists c > 0\, \forall \eps : c \le \omega_\eps,\\ 
                      &\forall p \ge 0\, \exists c_p > 0\, \forall\eps: \omega_\eps^p \le c_p\, \eps^{-1} \}.
\end{split}\eeq                      

\begin{definition}\label{def_ssc_symbols} 
 Let $m$ be a real number.
 The set of slow scale nets of symbols of order $m$ is defined by
\beq\label{def_ssc_symbols_mod}\begin{split}
 \Sscu^m(\Om\times\R^n) := \{ (a_\eps)_\eps\in {\mathcal{S}}^m[\Om\times\R^n]:\, &\forall K \Subset \Omega\ \exists (\omega_\eps)_\eps \in \Pi_\ssc\\
& \forall \alpha, \beta \in \N^n\, \exists c>0\,  \forall \eps:\, |a_\eps|^{(m)}_{K,\alpha,\beta} \le c\, \omega_\eps \},
\end{split}\eeq
 the negligible nets of symbols of order $m$ are the elements of 
\beq\label{def_negligible}\begin{split}
 \Nu^m(\Om\times\R^n) := \{  (a_\eps)_\eps\in {\mathcal{S}}^m[\Om\times\R^n]:\, &\forall K \Subset \Omega\, \forall \alpha, \beta \in \N^n\\
&\forall q\in\N, \exists c>0\,  \forall \eps:\, |a_\eps|^{(m)}_{K,\alpha,\beta} \le c\, \eps^q \}.
\end{split}
\eeq
The classes of the factor space 
\[
\Syscu^m(\Om\times\R^n) := \Sscu^m(\Om\times\R^n) / \Nu^m(\Om\times\R^n)
\]
are called \emph{slow scale generalized symbols of order $m$}.\\
Furthermore, let $\Nuinf(\Om\times\R^n):= \bigcap_m \Nu^m(\Om\times\R^n)$ be the negligible nets of order $-\infty$. The \emph{slow scale generalized symbols of refined order $m$} are given by 
\[
\Syscu^{m/-\infty}(\Om\times\R^n) := \Sscu^m(\Om\times\R^n) / \Nuinf(\Om\times\R^n).
\] 
\end{definition} 
Note that $\Syscu^{m/-\infty}(\Om\times\R^n)$ can be viewed as a finer partitioning of the classes in $\Syscu^m(\Om\times\R^n)$, in other words, if $a$ is a slow scale generalized symbol of order $m$ then $\forall (b_\eps)_\eps\in a$:
\begin{equation}
\label{kappa}
 \kappa((b_\eps)_\eps) := (b_\eps)_\eps +\Nuinf(\Om\times\R^n)\subseteq (b_\eps)_\eps + \Nu^{m}(\Om\times\R^n) = a.
\end{equation}
Slow scale generalized symbols enable us to design a particularly simple, yet sufficiently strong, notion of micro-ellipticity. 
\begin{definition}
\label{def_micro_ellipticity}
Let $a \in \Syscu^m(\Om\times\R^n)$ and $(x_0,\xi_0) \in \CO{\Om}$. We say that $a$ is slow scale micro-elliptic at $(x_0,\xi_0)$ if it has a representative $(a_\eps)_\eps$ satisfying the following: there is an relatively compact open neighborhood $U$ of $x_0$, a conic neighborhood $\Gamma$ of $\xi_0$, and $(r_\eps)_\eps , (s_\eps)_\eps$ in $\Pi_\ssc$ such that
\beq
\label{estimate_below}
 | a_\eps(x,\xi)| \ge \frac{1}{s_\eps} \lara{\xi}^m\qquad\qquad (x,\xi)\in U\times\Gamma,\, |\xi| \ge r_\eps,\, \eps \in (0,1].
\eeq
We denote by $\Ellsc(a)$ the set of all $(x_0,\xi_0) \in \CO{\Om}$ where $a$ is slow scale micro-elliptic.\\
If there exists $(a_\eps)_\eps \in\ a$ such that \eqref{estimate_below} holds at all points in $\CO{\Om}$ then the symbol $a$ is called \emph{slow scale elliptic}.
\end{definition}
Note that here the use of the attribute \emph{slow scale} refers to the appearance of the slow scale lower bound in \eqref{estimate_below}. This is a crucial difference to more general definitions of ellipticity given in \cite{Garetto:04, GGO:03, HOP:03}, whereas a similar condition was already used in \cite[Section 6]{HO:03} in a special case. In fact, due to the overall slow scale conditions in Definition \ref{def_ssc_symbols}, any symbol which is slow scale micro-elliptic at $(x_0,\xi_0)$ fulfills the stronger hypoellipticity estimates \cite[Definition 6.1]{GGO:03}; furthermore, \eqref{estimate_below} is stable under lower order (slow scale) perturbations. 
\begin{proposition}
\label{proposition_lower_order}
Let $(a_\eps)_\eps \in \Sscu^m(\Om\times\R^n)$ satisfy \eqref{estimate_below} in $U \times \Gamma \ni (x_0,\xi_0)$. Then
\begin{enumerate}
\item for all $\alpha,\beta \in \N^n$ there exists $(\lambda_\eps)_\eps \in \Pi_\ssc$ such that
\[
 |\partial^\alpha_\xi \partial^\beta_x a_\eps(x,\xi)| \le \lambda_\eps |a_\eps(x,\xi)| \lara{\xi}^{-|\alpha|}\qquad (x,\xi) \in U\times\Gamma,\, |\xi|\ge r_\eps,\, \eps\in(0,1];
\]
\item for all $(b_\eps)_\eps \in \Sscu^{m'}(\Om\times\R^n)$, $m'<m$, there exist $(r'_\eps)_\eps$, $(s'_\eps)_\eps \in \Pi_\ssc$ such that 
\[
|a_\eps(x,\xi)+b_\eps(x,\xi)| \ge \frac{1}{s'_\eps}\lara{\xi}^m\qquad (x,\xi) \in U\times\Gamma,\, |\xi|\ge r'_\eps,\, \eps\in(0,1].
\]
\end{enumerate}
\end{proposition}
\begin{proof}
Combining \eqref{estimate_below} with the seminorm estimates of $(a_\eps)_\eps \in \Sscu^m(\Om\times\R^n)$, we have that for $(x,\xi) \in U\times\Gamma$, $|\xi|\ge r_\eps$, $\eps\in(0,1]$
\[
 |\partial^\alpha_\xi \partial^\beta_x a_\eps(x,\xi)|\le c\, \omega_\eps \lara{\xi}^{m-|\alpha|}\le c\, \omega_\eps s_\eps |a_\eps(x,\xi)| \lara{\xi}^{-|\alpha|},
\]
so that assertion (i) holds with $\lambda_\eps = c\, \omega_\eps s_\eps$. To prove (ii), again by \eqref{estimate_below} for $(a_\eps)_\eps$ and the seminorm estimates for $(b_\eps)_\eps$ we obtain 
\[
|a_\eps(x,\xi)+b_\eps(x,\xi)|\ge \frac{1}{s_\eps}\lara{\xi}^m - c\, \omega_\eps \lara{\xi}^{m'}
 = \lara{\xi}^m \big( \frac{1}{s_\eps} - c\, \omega_\eps \lara{\xi}^{m'-m}),
\]
which is bounded from below by $\lara{\xi}^m / {2s_\eps}$ whenever $(x,\xi)\in U\times\Gamma$ with $|\xi| \ge \max( r_\eps , (2c\, \omega_\eps s_\eps)^{1/(m-m')})$.
\end{proof}
\begin{remark}\leavevmode
\label{remark_ellipticity}
\begin{trivlist}
\item[(i)] In case of classical symbols the notion of slow scale
micro-ellipticity coincides with the classical one, which equivalently is
defined as the set of noncharacteristic points. Indeed, if
$a\in\mathcal{S}^m(\Om\times\R^n)$ and $(a_\eps)_\eps$ is a representative of
the class of $a$ in $\Syscu^m(\Om\times\R^n)$ satisfying
\eqref{estimate_below} then for any $q\in\N$
\[
|a(x,\xi)| \ge |a_\eps(x,\xi)|-|(a-a_\eps)(x,\xi)|\ge \lara{\xi}^m(\frac{1}{s_\eps}-c\eps^q),
\]
where we are free to fix $\eps$ small enough such that the last factor is
bounded away from $0$. In particular, we have that $\compl{\Ellsc(a)} =
\Char(a(x,D))$. 
\item[(ii)] The same Definition \ref{def_micro_ellipticity} can be applied to
symbols of refined order.In that case, by Proposition
\ref{proposition_lower_order}, \eqref{estimate_below} will hold for any
representative once it is known to hold for one. Moreover, if
$a\in\Syscu^{m/-\infty}(\Om\times\R^n)$ and $\compl{\Ellsc(a)} = \emptyset$
then $a$ is slow scale elliptic. 
\end{trivlist}
\end{remark}
Thanks to the previous proposition the simple slow scale ellipticity condition in Definition \ref{def_micro_ellipticity} already guarantees the existence of a parametrix. For the proof we refer to \cite[Section 6]{GGO:03}; note that an inspection of the construction shows that the symbol of the parametrix has uniform growth $\eps^{-1}$ over all compact sets. Regular symbols are introduced in Remark \ref{rem_special}.
\begin{theorem}
\label{theorem_parametrix}
Let $a$ be a slow scale elliptic symbol of order $m$. Then there exists a properly supported pseudodifferential operator with regular symbol $p \in \Syru^{-m}(\Om\times\R^n)$ such that for all $u\in\Gc(\Om)$
\[
\begin{split}
 a(x,D) \circ p(x,D)u &= u + Ru,\\
 p(x,D) \circ a(x,D)u &= u + Su,
\end{split}
\]
where $R$ and $S$ are operators with regular kernel.
\end{theorem}
Note that if $a(x,D)$ is properly supported then the operators $R$ and $S$ are properly supported too and the previous equalities are valid for all $u$ in $\G(\Om)$. In this situation, combining the construction of a parametrix with the pseudolocality property (see the Appendix for details), we obtain that $\singsupp_g(a(x,D)u) = \singsupp_g(u)$ for all $u\in\G(\Om)$.\\ 
In the sequel $\Oprop{m}(\Om)$ denotes the set of all properly supported operators $a(x,D)$ where $a$ belongs to $\Syscu^m(\Om\times\R^n)$. We are now in a position to introduce a way to measure regularity of Colombeau generalized functions mimicking the original definition of the distributional wave front set in \cite{Hoermander:71} based on characteristic sets. As a matter of fact, the set constructed below as the complement of the slow scale micro-ellipticity regions will turn out to be the generalized wave front set in the sense of \cite{DPS:98, Hoermann:99}.
\begin{definition}
Let $u\in\G(\Om)$. We define
\beq
\label{def_W_ssc}
 \Wsc(u) :=\hskip-5pt \bigcap_{\substack{a(x,D)\in\,\Oprop{0}(\Om)\\[0.1cm] a(x,D)u\, \in\, \Ginf(\Om)}}\hskip-5pt \compl{\Ellsc(a)}.
\eeq
\end{definition}
\begin{remark}
\label{remark_any_order}
Note that the standard procedure of lifting symbol orders with $(1-\Delta)^{m/2}$ easily shows that we may as well take the intersection over operators $a(x,D)\in\Oprop{m}(\Om)$ in \eqref{def_W_ssc}. The same holds for similar constructions introduced throughout the paper.
\end{remark}
Since $\compl{\Ellsc(a)}$ is a closed conic set, $\Wsc(u)$ is a closed conic subset of $\CO{\Om}$ as well. Moreover, recalling that given $v\in\Ginf(\Om)$ and $a(x,D)$ properly supported, $a(x,D)(u+v)\in\Ginf(\Om)$ iff $a(x,D)u\in\Ginf(\Om)$, we have $\Wsc(u+v)=\Wsc(u)$.\\
We present now a first alternative way to define $\Wsc(u)$, which will be useful in course of our exposition. Denote by $\Oprop{m/-\infty}(\Om)$ the set of all properly supported operators $a(x,D)$ where $a\in\Syscu^{m/-\infty}(\Om\times\R^n)$; one can prove that
\beq
\label{alt_def_W_ssc} 
 \Wsc(u) =\hskip-5pt \bigcap_{\substack{a(x,D)\in\,\Oprop{0/-\infty}(\Om)\\[0.1cm] a(x,D)u\, \in\, \Ginf(\Om)}}\hskip-5pt \compl{\Ellsc(a)}.
\eeq
In fact, the crucial point is to observe that we do not change a pseudodifferential operator with generalized symbol $a$ by adding negligible nets of symbols of the same order or of order $-\infty$. To be more precise, if $a(x,D)\in\Oprop{0}(\Om)$, $(x_0,\xi_0)\in\Ellsc(a)$ with $(a_\eps)_\eps$ satisfying \eqref{estimate_below}, then $b(x,\xi):=(a_\eps)_\eps +\Nu^{-\infty}(\Om\times\R^n)$ is slow scale micro-elliptic at $(x_0,\xi_0)$,  $b(x,D)\in\Oprop{0/-\infty}(\Om)$  since $b(x,D)\equiv a(x,D)$ and consequently we have the inclusion $\supseteq$ in \eqref{alt_def_W_ssc}. For the reverse inclusion if $a(x,D)\in\Oprop{0/-\infty}(\Om)$ with $(x_0,\xi_0)\in{\Ellsc(a)}$, it is clear that $b=(a_\eps)_\eps +\Nu^0(\Om\times\R^n)$ is a well-defined element of $\Syscu^0(\Om\times\R^n)$ and slow scale elliptic at $(x_0,\xi_0)$. Arguing as before we obtain \eqref{alt_def_W_ssc}.  
From a technical point of view, the most interesting aspect of \eqref{alt_def_W_ssc} is the stability of micro-ellipticity estimates under variations of the representatives of $a$, valid for symbols of refined order due to Proposition \ref{proposition_lower_order} .
\begin{proposition}
\label{prop_sing_supp}
Let $\pi:\CO{\Om}\to\Om:(x,\xi)\to x$. For any $u\in\G(\Om)$, 
\[
\pi(\Wsc(u)) = \singsupp_\g(u).
\]
\end{proposition}
\begin{proof}
We first prove that $\Om\setminus\singsupp_\g(u) \subseteq \Om\setminus
\pi(\Wsc(u))$. Let $x_0\in\Om\setminus\singsupp_\g(u)$. There exists 
$\phi\in\Cinf_{\rm{c}}(\Om)$ such that $\phi(x_0)=1$ and $\phi u\in\Ginf(\Om)$. The multiplication operator $\phi(x,D):\Gc(\Om)\to\Gc(\Om):u\to\phi u$ is properly supported with symbol $\phi\in S^0(\Om\times\R^n)\subseteq\Syscu^0(\Om\times\R^n)$, which is micro-elliptic (slow scale micro-elliptic) at $(x_0,\xi_0)$ for each $\xi_0\neq 0$. Therefore, for all $\xi_0\neq 0$ we have that $(x_0,\xi_0)\in\compl{\Wsc(u)}$ i.e. $x_0\in\Om\setminus \pi(\Wsc(u))$.\\
To show the opposite inclusion, let $x_0\in\Om\setminus\pi(\Wsc(u))$. Then for all $\xi\neq 0$ there exists $a\in\Syscu^{0/-\infty}(\Om\times\R^n)$ slow scale micro-elliptic at $(x_0,\xi)$ such that $a(x,D)$ is properly supported and $a(x,D)u\in\Ginf(\Om)$. Since $S_{x_0} := \{x_0\}\times \{\xi: |\xi|=1\}$ is a compact subset of $\Om\times\R^n$, there exist $a_1,...,a_N\in\Syscu^{0/-\infty}(\Om\times\R^n)$, $U$ a relatively compact open neighborhood of $x_0$ and $\Gamma_i$ conic neighborhoods of $\xi_i$ with $|\xi_i|=1$ ($i=1,\ldots,N$), with the following properties: $a_i$ is slow scale micro-elliptic in $U\times\Gamma_i$, $S_{x_0}\subseteq U\times\cup_{i=1}^N\Gamma_i$, $a_i(x,D)$ is properly supported, and $a_i(x,D)u\in\Ginf(\Om)$. Consider the properly supported pseudodifferential operator $A:=\sum_{i=1}^N a_i(x,D)^\ast a_i(x,D)$. By Theorem \ref{theorem_calculus} in the Appendix we may write $A=\sigma(x,D)\in\Oprop{0/-\infty}(\Om)$, and combining assertions $(ii)$ and $(iii)$ of the same theorem, we have that 
\beq
\label{differ_sigma}
 \sigma-\sum_{i=1}^N|a_i|^2 \in \Syscu^{-1/-\infty}(\Om\times\R^n).
\eeq
Since $a_i(x,D)u\in\Ginf(\Om)$ and each $a_i(x,D)^\ast$ maps $\Ginf(\Om)$ into $\Ginf(\Om)$, we conclude that $\sigma(x,D)u\in\Ginf(\Om)$.\\
It is clear that $\sum_{i=1}^N|a_i|^2$ is slow scale elliptic in $U$. In fact, every $\xi\neq 0$ belongs to some $\Gamma_i$ and, given $(s_{i,\eps})_\eps, (r_{i,\eps})_\eps$ satisfying \eqref{estimate_below} for $(a_{i,\eps})_\eps$, $s_\eps:=\max_{i}(s^2_{i,\eps})$, $r_\eps:=\max(r_{i,\eps})$, we get
\[
\sum_{i=1}^N|a_{i,\eps}(x,\xi)|^2\ge \frac{1}{s_\eps}\qquad\qquad x\in U,\, |\xi|\ge r_\eps,\, \eps\in(0,1].
\]
Let $U'\subset U''\subset U$ be open neighborhoods of $x_0$, $\overline{U'}\subset U^{''}$, $\overline{U''}\subset U$ and $\phi\in\Cinf(\Om)$, $0\le\phi\le 1$ such that $\phi=0$ on $U'$ and $\phi=1$ on $\Om\setminus U''$. By construction $b(x,\xi):=\phi(x)+\sigma(x,\xi)\in\Syscu^{0/-\infty}(\Om\times\R^n)$, $b(x,D)\in\Oprop{0/-\infty}(\Om)$ and $b(x,D)u_{\vert_{U'}} = \phi u_{\vert_{U'}} + \sigma(x,D)u_{\vert_{U'}} = \sigma(x,D)u_{\vert_{U'}} \in \Ginf(U')$.  Since $\sum_{i=1}^N|a_i|^2$ is slow scale elliptic in $U$, with a positive real valued representative, and $\phi$ is identically $1$ outside $U''$, we have that $\phi+\sum_{i=1}^N|a_i|^2$ is slow scale elliptic in $\Om$. By \eqref{differ_sigma}, and application of Proposition \ref{proposition_lower_order}(ii), $b$ itself is slow scale elliptic in $\Om$. Then, by using a parametrix for $b(x,D)$ we conclude that $\singsupp_g(b(x,D)u) = \singsupp_g(u)$ and consequently $U'\cap\singsupp_g(u) =\emptyset$, which completes the proof.  
\end{proof}

\section{Pseudodifferential characterization of the\\ generalized wave front set}
This section is devoted to the proof that $\Wsc(u)$ coincides with the generalized wave front set of $u$. Our approach will follow the lines of reasoning in \cite[Chapter 8]{Folland:95} and \cite[Section 18.1]{Hoermander:V3}. The main tool will be a generalization of the micro-support of a regular generalized symbol of order $m$ or refined order $m$.
\begin{definition}
\label{def_micro_support}
Let $a\in\Syru^{m}(\Om\times\R^n)$ and $(x_0,\xi_0)\in\CO{\Om}$. The symbol $a$ is $\Ginf$-smoothing at $(x_0,\xi_0)$ if there exist a representative $(a_\eps)_\eps$ of $a$, a relatively compact open neighborhood $U$ of $x_0$, a conic neighborhood $\Gamma$ of $\xi_0$, and a natural number $N$ such that
\beq
\label{estimate_micro}
\begin{split}
\forall m\in \R\  \forall\alpha ,\beta\in\mathbb{N}^n\ &\exists c>0\ \forall(x,\xi)\in U\times\Gamma\  \forall\eps\in(0,1]:\\
&|\partial^\alpha_\xi\partial^\beta_x a_\eps(x,\xi)|\le c\lara{\xi}^m \eps^{-N}.
\end{split}
\eeq
We define the \emph{generalized microsupport} of $a$, denoted by $\mu\supp_g(a)$, as the complement of the set of points $(x_0,\xi_0)$ where $a$ is $\Ginf$-smoothing.\\
If $a\in\Syru^{m/-\infty}(\Om\times\R^n)$ then we denote by $\mufty(a)$ the complement of the set points $(x_0,\xi_0)\in\CO{\Om}$ where \eqref{estimate_micro} holds for some representative of $a$.
\end{definition}
\begin{remark}\leavevmode
\label{remark_micro_support}
\begin{trivlist}
\item[(i)] Any $a\in\Syru^{-\infty}(\Om\times\R^n)$, i.e., a regular generalized symbol of order $-\infty$, has empty generalized microsupport.
\item[(ii)] In the case of $a\in\Syru^{m/-\infty}(\Om\times\R^n)$, every representative is of the form $(a_\eps)_\eps +(n_\eps)_\eps$, where $(a_\eps)_\eps\in\Sru^m(\Om\times\R^n)$ and $(n_\eps)_\eps\in\Nu^{-\infty}(\Om\times\R^n)$, and \eqref{estimate_micro} holds for any representative once it is known to hold for one. As a consequence, if $a$ is a classical symbol of order $m$ considered as an element of $\Syru^{m/-\infty}(\Om\times\R^n)$, then its generalized  microsupport $\mufty(a)$ equals the classical one.
\item[(iii)] If $a\in\Syru^{m/-\infty}(\Om\times\R^n)$ and $\mu\supp_\g(a)=\emptyset$ then $a\in\Syru^{-\infty}(\Om\times\R^n)$. In fact, chosen any $(a_\eps)_\eps\in a$, for all $x_0\in\Om$ the compact set $S_{x_0}$ can be covered by $U\times\cup_{i=1}^N\Gamma_i$, with $U$ and $\Gamma_i$ such that \eqref{estimate_micro} is valid. Hence there exists $N\in\mathbb{N}$ such that for all orders $m\in\mathbb{R}$ and for all $\alpha,\beta$
\[
|\partial^\alpha_\xi\partial^\beta_x a_\eps(x,\xi)|\le c\lara{\xi}^{m-|\alpha|}\eps^{-N}\qquad\qquad (x,\xi)\in U\times\R^n,\ \eps\in(0,1].
\]
\item[(iv)] If $a\in\Syru^m(\Om\times\R^n)$ and $\kappa$ is the quotient map $\Sru^{m}(\Om\times\R^n)$ onto $\Syru^{m/-\infty}(\Om\times\R^n)$ then
\beq
\label{intmusupp}
 \mu\supp_\g(a) = \bigcap_{(a_\eps)_\eps \in a} \mufty\big(\,\kappa\big((a_\eps)_\eps\big)\,\big).
\eeq 
 Observe first that for all $(a_\eps)_\eps\in a$ we have $\mu\supp_g(a) \subseteq \mufty(\kappa((a_\eps)_\eps))$.
On the other hand, if $(x_0,\xi_0)\not\in  \mu\supp_g(a)$ then $(x_0,\xi_0)\not\in \mufty(\kappa((a_\eps)_\eps))$ for some $(a_\eps)_\eps\in a$.
\end{trivlist}
\end{remark}
Similarly, as in the previous section, we introduce the notations $\Opropr{m}(\Om)$ and $\Opropr{m/-\infty}(\Om)$ for the sets of all properly supported operators $a(x,D)$ with symbol in $\Syru^{m}(\Om\times\R^n)$ and $\Syru^{m/-\infty}(\Om\times\R^n)$ respectively.
\begin{proposition}
\label{micro_support_product}
Let $a(x,D)\in\Opropr{m/-\infty}(\Om)$ and $b(x,D)\in\Opropr{m'/-\infty}(\Om)$. Then, there exists a symbol $a\sharp b\in\Syru^{m+m'/-\infty}(\Om\times\R^n)$ such that $a(x,D)\circ b(x,D) = a\sharp b(x,D)\in\Opropr{m+m'/-\infty}(\Om)$ and
\beq
\label{micro_product}
 \mufty(a\sharp b) \subseteq \mufty(a) \cap \mufty(b).
\eeq
In the same situation, without regarding refined orders we have that
\beq
\label{micro_product_2}
\mu\supp_g(a\sharp b) \subseteq \mu\supp_g(a) \cap \mu\supp_g(b).
\eeq
\end{proposition}
\begin{proof}
Theorem \ref{theorem_calculus} provides the existence of $a\sharp b\in\Syru^{m+m'/-\infty}(\Om\times\R^n)$ such that $a(x,D)\circ b(x,D) = a\sharp b(x,D)\in\Opropr{m+m'/-\infty}(\Om)$ with asymptotic expansion $a\sharp b \sim \sum_\gamma \partial^\gamma_\xi a D^\gamma_x b/\gamma!$. Assume that $(x_0,\xi_0)\not\in{\mufty(a)}$. By \eqref{estimate_micro} and the symbol properties of $(b_\eps)_\eps\in b$, we obtain the following estimate valid on some neighborhood $U\times\Gamma$ of $(x_0,\xi_0)$:
\beq
\label{estimate1}
\exists N\in\N\ \forall m\in\R\ \forall\alpha,\beta\in\N^n\qquad |\partial^\alpha_\xi\partial^\beta_x(\partial^\gamma_\xi a_\eps D^\gamma_x b_\eps)(x,\xi)|\le c\lara{\xi}^m\eps^{-N}.
\eeq
Since $a\sharp b\sim\sum_\gamma\partial^\gamma_\xi a D^\gamma_x b/\gamma !$ we have that for any $(d_\eps)_\eps\in a\sharp b$ and $h\ge 1$, the difference $((d_\eps)_\eps -\sum_{|\gamma|\le h-1}\frac{1}{\gamma !}\partial^\gamma_\xi a_\eps D^\gamma_x b_\eps)_\eps$ is an element of $\Sru^{m+m'-h}(\Om\times\R^n)$ and of growth type $\eps^{-M}$ on $\overline{U}$, for some $M\in\N$ independent of $h$. This together with \eqref{estimate1} implies that $(x_0,\xi_0)$ does not belong to ${\mufty(a\sharp b)}$. The proof of relation \eqref{micro_product_2} is similar.
\end{proof}
We recall a technical lemma proved in \cite[Proposition (8.52)]{Folland:95} which will be useful in the sequel.
\begin{lemma}
\label{lemma_folland}
Suppose $(x_0,\xi_0)\in\CO{\Om}$, $U$ is a relatively compact open neighborhood of $x_0$, $\Gamma$ is a conic neighborhood of $\xi_0$. There exists $p\in{S}^0(\Om\times\R^n)$ such that $0\le p\le 1$, $\supp(p)\subseteq U\times\Gamma$ and $p(x,\xi)=1$ if $(x,\xi)\in U'\times\Gamma'$ and $|\xi|\ge 1$, where $U'\times\Gamma'$ is a smaller conic neighborhood of $(x_0,\xi_0)$. In particular, $p$ is micro-elliptic at $(x_0,\xi_0)$ and $\mu\supp(p)\subseteq U\times\Gamma$.
\end{lemma}
\begin{remark}
\label{remark_folland}
The proof of the above lemma actually shows that for each conic neighborhood $\Gamma$ of $\xi_0$ there exists $\tau(\xi)\in{S}^0(\Om\times\R^n)$ such that $0\le \tau\le 1$, $\supp(\tau)\subseteq\Gamma$ and $\tau(\xi)=1$ in some conic neighborhood $\Gamma'$ of $\xi_0$ when $|\xi|\ge 1$.\\ 
Note that after multiplying $p(x,\xi)$ in Lemma \ref{lemma_folland} with a proper cut-off function, we obtain a properly supported operator $q(x,D)\in\Psi^{0}(\Om)$ whose symbol is micro-elliptic at $(x_0,\xi_0)$ and $\mu\supp(q)=\mu\supp(p)$.
\end{remark}
\begin{theorem}
\label{theorem_micro_support}
For any $a(x,D)\in\Opropr{m}(\Om)$ and $u\in\G(\Om)$
\beq
\label{inclusion}
 \Wsc(a(x,D)u) \subseteq \Wsc(u) \cap \mu\supp_\g(a).
\eeq
Similarly, if $a(x,D)\in\Opropr{m/-\infty}(\Om)$ then 
\beq
\label{inclusion_1}
 \Wsc(a(x,D)u) \subseteq \Wsc(u) \cap \mufty(a).
\eeq
\end{theorem}
\begin{proof}
We first prove the assertion \eqref{inclusion_1} in two steps.

\emph{Step 1:}\enspace\enspace ${\Wsc(a(x,D)u)}\subseteq \mufty(a)$.

If $(x_0,\xi_0)\not\in\mufty(a)$ then \eqref{estimate_micro} holds on some $U\times\Gamma$, and by Lemma \ref{lemma_folland} we find $q\in{S}^{0}(\Om\times\R^n)\subseteq\Syru^{0/-\infty}(\Om\times\R^n)$, which is micro-elliptic at $(x_0,\xi_0)$ with $\mu\supp(q)\subseteq U\times\Gamma$. Apply Proposition \ref{micro_support_product} to obtain $q(x,D)a(x,D)= q\sharp a(x,D)\in\Opropr{m/-\infty}(\Om)$ and $\mufty(q\sharp a)\subseteq \mufty(q)\cap \mufty(a)\subseteq (U\times\Gamma) \cap \mufty(a) = \emptyset$. Remark \ref{remark_micro_support}(iii) shows that $q\sharp a \in\Syru^{-\infty}(\Om\times\R^n)$ and therefore $q(x,D)a(x,D)$ is a properly supported pseudodifferential operator with regular kernel. This implies $q(x,D)(a(x,D)u)\in\Ginf(\Om)$ and hence $(x_0,\xi_0)\not\in{\Wsc(a(x,D)u)}$.

\emph{Step 2:}\enspace\enspace ${\Wsc(a(x,D)u)}\subseteq \Wsc(u)$.

Let $(x_0,\xi_0)\not\in{\Wsc(u)}$ then by \eqref{alt_def_W_ssc} there exists $p(x,D)\in\Oprop{0/-\infty}(\Om)$ such that $p$ is slow scale micro-elliptic at $(x_0,\xi_0)$ and $p(x,D)u\in\Ginf(\Om)$.

\emph{Claim}: There exists $r(x,D)\in\Oprop{0/-\infty}(\Om)$ such that $r$ is micro-elliptic at $(x_0,\xi_0)$ and there exists $s(x,D)\in\Opropr{m/-\infty}(\Om)$ such that $r(x,D)a(x,D)u-s(x,D)p(x,D)^\ast p(x,D)u$ belongs to $\Ginf(\Om)$.

Assuming for the moment that the claim is proved, we show that it completes the proof of the theorem. The operators $s(x,D)$ and $p(x,D)^\ast$ map $\Ginf(\Om)$ into itself, hence we obtain that $s(x,D)p(x,D)^\ast p(x,D)u\in\Ginf(\Om)$. The claim implies that $r(x,D)a(x,D)u\in\Ginf(\Om)$ and $(x_0,\xi_0)\not\in{\Wsc(a(x,D)u)}$, since $r$ is micro-elliptic at $(x_0,\xi_0)$.\\

To prove the claim we construct a slow scale elliptic symbol based on $p$. Let $(p_\eps)_\eps$ be a representative of $p$ satisfying \eqref{estimate_below} in a conic neighborhood $U\times\Gamma$ of $(x_0,\xi_0)$. By Lemma \ref{lemma_folland} there is $\psi\in{S}^0(\Om\times\R^n)$, $0\le\psi\le 1$, with $\supp(\psi)\subseteq U\times\Gamma$ and identically $1$ in a smaller conic neighborhood $U'\times\Gamma'$ of $(x_0,\xi_0)$ if $|\xi|\ge 1$. The net $(1+|p_\eps|^2-\psi)_\eps$ belongs to $\Sscu^0(\Om\times\R^n)$, and by construction of $\psi$, satisfies \eqref{estimate_below} at all points in $\CO{\Om}$. The pseudodifferential calculus for slow scale generalized symbols of refined order guarantees the existence of $\sigma\in\Syscu^{0/-\infty}(\Om\times\R^n)$ such that $p(x,D)^\ast p(x,D) = \sigma(x,D)\in\Oprop{0/-\infty}(\Om)$ and $\sigma-|p|^2\in\Syscu^{-1/-\infty}(\Om\times\R^n)$. Application of Proposition \ref{proposition_lower_order}(ii) to this situation yields that the symbol $1+\sigma-\psi\in\Syscu^{0/-\infty}(\Om\times\R^n)$ is slow scale elliptic and coincides with $\sigma$ in a conic neighborhood $U'\times\Gamma'$ of $(x_0,\xi_0)$ for $|\xi|\ge 1$.\\
Take a proper cut-off function $\chi$ and define a pseudodifferential operator via the slow scale amplitude $\chi(x,y)(1+\sigma(x,\xi)-\psi(x,\xi))$. By Theorem \ref{theorem_calculus} it can be written in the form  $b(x,D)\in\Oprop{0/-\infty}(\Om)$, where $b-(1+\sigma-\psi)\in\Syscu^{-\infty}(\Om\times\R^n)$. In other words, $b$ itself is slow scale elliptic and $\mufty(b) = \mufty(1+\sigma-\psi)$. In particular, since $\sigma$ and $1+\sigma-\psi$ coincide on $U'\times\Gamma'$, $|\xi|\ge 1$, we have that $\mufty(b-\sigma)\cap ( U'\times\Gamma') = \emptyset$. Theorem \ref{theorem_parametrix} gives a parametrix $t(x,D)\in\Opropr{0/-\infty}(\Om)$ for $b(x,D)$, i.e., the operators $b(x,D)t(x,D)-I$ and $t(x,D)b(x,D)-I$ have regular kernel. Let $r(x,D)\in\Oprop{0/-\infty}(\Om)$ be an operator constructed as in Lemma \ref{lemma_folland}, with classical symbol $r$, micro-elliptic at $(x_0,\xi_0)$, and $\mu\supp(r)\subseteq U'\times\Gamma'$.\\
We show that $r(x,D)$ and $s(x,D):=r\sharp(a\sharp t)(x,D)=r(x,D)a(x,D)t(x,D)\in\Opropr{m/-\infty}(\Om)$ satisfy the assertions of the claim. We rewrite the difference $r(x,D)a(x,D)u-s(x,D)p(x,D)^\ast p(x,D)u$ as 
\begin{multline*}
r(x,D)a(x,D)\big(u-t(x,D)b(x,D)u\big)\\
 + r(x,D)a(x,D)t(x,D)\big(b(x,D)-\sigma(x,D)\big)u.
\end{multline*}
Here, the first summand is in $\Ginf(\Om)$ due to the fact that $t(x,D)$ is a parametrix for $b(x,D)$ and the mapping properties of $r(x,D)a(x,D)$. An iterated application of Proposition \ref{micro_support_product} to the second summand shows that it can be written with a regular symbol of refined order $m$, having generalized micro-support contained in the region $\mu\supp(r)\cap\mufty(b-\sigma)\subseteq (U'\times\Gamma') \cap \mufty(b-\sigma) =\emptyset$. Hence it has  smoothing generalized symbol and therefore the claim is proven.\\

Finally we prove the assertion \eqref{inclusion}. Let $(a_\eps)_\eps\in a \in \Syru^{m}(\Om\times\R^n)$ and consider the corresponding symbol $\kappa((a_\eps)_\eps)=(a_\eps)_\eps +\Nu^{-\infty}(\Om\times\R^n)\in\Syru^{m/-\infty}(\Om\times\R^n)$. Then $\kappa((a_\eps)_\eps)(x,D)\in\Opropr{m/-\infty}(\Om)$ and $\kappa((a_\eps))(x,D)=a(x,D)$. Theorem \ref{theorem_micro_support} applied to all $(a_\eps)_\eps \in a$ yields 
\[
\bigcap_{(a_\eps)_\eps \in a}\hskip-4pt \Wsc(\kappa((a_\eps)_\eps)(x,D)u)\ \subseteq\ \Wsc(u)\ \cap\ \bigcap_{(a_\eps)_\eps \in a}\hskip-4pt \mufty(\kappa((a_\eps)_\eps))
\]
i.e.,
\[
\Wsc(a(x,D)u)\ \subseteq\ \Wsc(u)\ \cap\ \bigcap_{(a_\eps)_\eps \in a}\hskip-4pt \mufty(\kappa((a_\eps)_\eps)),
\]
which completes the proof by Remark \ref{remark_micro_support}(iv).  
\end{proof}
\begin{corollary}
\label{corollary_2}
Let $a(x,D)\in\Oprop{m}(\Om)$ where $a$ is a slow scale elliptic symbol. Then for any $u\in\G(\Om)$
\beq
\label{equality_corollary_ell}
 \Wsc(a(x,D)u) = \Wsc(u).
\eeq
\end{corollary}
\begin{proof}
Since $a\in\Syscu^m(\Om\times\R^n)\subseteq\Syru^m(\Om\times\R^n)$, Theorem \ref{theorem_micro_support} implies that $\Wsc(a(x,D)u)\subseteq \Wsc(u)$. Let $p(x,D)$ be a parametrix for $a(x,D)$ as in Theorem \ref{theorem_parametrix} then $u = p(x,D)a(x,D)u+v$, where $v\in\Ginf(\Om)$. Therefore $\Wsc(u)=\Wsc(p(x,D)a(x,D)u)$ and Theorem \ref{theorem_micro_support} applied to $p(x,D)\in\Opropr{-m}(\Om)$ gives $\Wsc(u)\subseteq\Wsc(a(x,D)u)$.
\end{proof}
The statements of the above theorem and corollary are valid for pseudodifferential operators which are not necessarily properly supported, if we consider instead compactly supported generalized functions.\\

As in the classical theory (cf.\ \cite{Folland:95}), we introduce a
notion of microsupport for operators. However, in the case of generalized
pseudodifferential operators we are cautious to distinguish the corresponding
notions for symbols and operators by a slight change in notation, thereby
taking into account the non-injectivity when mapping symbols to operators
(cf.\ \cite{GGO:03}).
\begin{definition}
\label{def_micro_supp_op}
Let $A$ be any properly supported pseudodifferential operator with regular symbol. We define the \emph{generalized microsupport} of $A$ by 
\beq
\label{micro_supp_op}
\Bigmu\supp_g(A)\ :=\ \bigcap_{\substack{ a(x,D)\in\Opropr{m}(\Om),\\[0.1cm]\hskip-10pt a(x,D)=A}}\mu\supp_g(a).
\eeq
\end{definition}

Now we have all the technical tools at hand which enable us to identify $\Wsc(u)$ as the generalized wave front set $\WFg(u)$. 
\begin{theorem}
\label{theorem_wave_front}
For all $u\in\G(\Om)$
$$\Wsc(u)=\WFg(u).$$  
\end{theorem}
\begin{proof}
By definition of $\WFg(u)$ the assertion which we are going to prove is the following: $(x_0,\xi_0)\not\in{\Wsc(u)}$ iff there exists a representative $(u_\eps)_\eps$ of $u$, a cut-off function $\phi\in\mathcal{C}^\infty_{\rm{c}}(\Om)$ with $\phi(x_0)=1$, a conic neighborhood $\Gamma$ of $\xi_0$ and a number $N$ such that for all $l\in\mathbb{R}$
\beq
\label{estimate_wave_front}
\sup_{\xi\in\Gamma}\,\lara{\xi}^l |\widehat{\phi u_\eps}(\xi)| = O(\eps^{-N})\qquad \text{as}\ \eps\to 0.
\eeq
We first show sufficiency, that is if $(x_0,\xi_0)\in\CO{\Om}$ satisfies \eqref{estimate_wave_front} then it does not belong to $\Wsc(u)$. As noted in Remark \ref{remark_folland}, there exists $p(\xi)\in{S}^0(\Om\times\R^n)$ with $\supp(p)\subseteq\Gamma$, which is identically $1$ in a conical neighborhood $\Gamma'$ of $\xi_0$ when $|\xi|\ge 1$. We recall that taking a typical proper cut-off $\chi$, by Theorem \ref{theorem_pre_calculus}, we can write the properly supported pseudodifferential operator with amplitude $\chi(x,y)p(\xi)\phi(y)$ in the form $\sigma(x,D)\in\Oprop{0/-\infty}(\Om)$, where $\sigma(x,\xi)-p(\xi)\phi(x)\in{S}^{-1}(\Om\times\R^n)$; in particular, $\sigma(x,D)v-p(D)(\phi v)\in\Ginf(\Om)$ for all $v\in\G(\Om)$. By assumption, $p(\xi)\phi(x)$ is micro-elliptic at $(x_0,\xi_0)$, the symbol $\sigma$ is micro-elliptic there and from \eqref{estimate_wave_front} we have that 
\[
p(D)\phi(x,D)u = \biggl[\biggl(\int_{\R^n}e^{ix\xi}p(\xi)\widehat{\phi u_\eps}(\xi)\, \dslash\xi\biggr)_\eps\biggr] \in \Ginf(\Om).
\]
Since $\sigma(x,D)u-p(D)(\phi u)\in\Ginf(\Om)$ we obtain that $\sigma(x,D)u\in\Ginf(\Om)$ and hence $(x_0,\xi_0)\not\in{\Wsc(u)}$.\\
Conversely, suppose $(x_0,\xi_0)\not\in{\Wsc(u)}$. There is an open neighborhood $U$ of $x_0$ such that $(x,\xi_0)\in\compl{\Wsc(u)}$ for all $x\in U$. Choose $\phi\in\Cinf_{\rm{c}}(U)$ with $\phi(x_0)=1$ and define 
\[
\Sigma :=\{ \xi\in\R^n\setminus 0:\ \ \exists x\in\Om\ (x,\xi)\in\Wsc(\phi u)\}.
\]
From Theorem \ref{theorem_micro_support} we have that $\Wsc(\phi u)\subseteq \Wsc(u) \cap (\supp(\phi)\times\R^n\setminus 0)$ and therefore $\xi_0\notin\Sigma$. Moreover, since $\Wsc(\phi u)$ is closed and  conic, the $\Sigma$ itself is a closed conic subset of $\R^n\setminus 0$. Again, by Remark \ref{remark_folland} there is a symbol $p(\xi)\in{S}^0(\Om\times\R^n)$ such that $0\le p\le 1$, $p(\xi)=1$ in a conic neighborhood $\Gamma$ of $\xi_0$ when $|\xi|\ge 1$ and $p(\xi)=0$ in a conic neighborhood $\Sigma_0$ of $\Sigma$. By construction $\mu\supp(p) \cap (\Om \times {\Sigma_0}) = \emptyset$ and $\Wsc(\phi u)\subseteq\Omega\times\Sigma$. Therefore, $\Wsc(p(D)\phi u)\subseteq \Wsc(\phi u)\cap \mu\supp(p)=\emptyset$ and by Proposition \ref{prop_sing_supp} we conclude that $p(D)\phi u\in\Ginf(\Om)$. In terms of representatives $(u_\eps)_\eps \in u$, this means that  
\beq
\label{representative}
\biggl(\int_{\R^n}e^{ix\xi}p(\xi)\widehat{\phi u_\eps}(\xi)\, \dslash\xi = (\phi u_\eps \ast {p}^{\vee})\biggr)_\eps \in\EMinf.
\eeq
Note that ${p}^\vee$ is a Schwartz function outside the origin, i.e., we have for all $\delta>0$ and $\alpha,\beta\in\N^n$ that $\sup_{|x|>\delta}|x^\alpha\partial^\beta {p}^\vee(x)|<\infty$ (\cite[Theorem (8.8a)]{Folland:95}). If ${\rm{dist}}(x,\supp(\phi))>\delta$ this yields $\partial^\alpha(\phi u_\eps\ast{p}^\vee)(x)=\int_{|y|>\delta}\phi u_\eps(x-y)\partial^\alpha{p}^\vee(y)\, dy$ and for all $l>0$ the estimate
\beq
\label{estimate_charac}
\begin{split}
&\lara{x}^l|\partial^\alpha(\phi u_\eps\ast p^\vee)(x)|\\
&\le c_1\int_{|y|>\delta}|\phi u_\eps(x-y)|\lara{x-y}^l\lara{y}^l|\partial^\alpha{p}^\vee(y)|\, dy\\
&\le c_2 \sup_{z\in\supp(\phi)}|u_\eps(z)|,
\end{split}
\eeq
uniformly for such $x$. When $\delta$ is chosen small enough, the set of points $x$ with ${\rm{dist}}(x,\supp(\phi))\le\delta$ is compact in $\Om$, and from \eqref{representative} we have that there exists $M\in\N$ such that for all $l\in\R$, $\sup \lara{x}^l|\partial^\alpha(\phi u_\eps\ast{p}^\vee)(x)| = O(\eps^{-M})$ as $\eps\to 0$, where the supremum is taken over $\{ x: {\rm{dist}}(x,\supp(\phi))\le\delta \}$. Hence we have shown that $(\phi u_\eps\ast{p}^\vee)_\eps\in\ESinf(\R^n)$. The Fourier transform maps $\ESinf(\R^n)$ into $\ESinf(\R^n)$, which implies that $(p(\xi)\widehat{\phi u_\eps}(\xi))_\eps\in\ESinf(\R^n)$. Since $p(\xi)=1$ in a conic neighborhood of $\xi_0$, the proof is complete.
\end{proof}
Combing Theorem \ref{theorem_micro_support} with Theorem \ref{theorem_wave_front} the following corollary is immediate from Definition \ref{def_micro_supp_op}.
\begin{corollary}
\label{corollary_wave_front}
For any properly supported operator $A$ with generalized regular symbol and $u\in\G(\Om)$ we have
\beq
\label{wave_front_op}
\WFg(Au) \subseteq \WFg(u) \cap \Bigmu\supp_\g(A).
\eeq
\end{corollary}

While the purpose of the foregoing discussion was to prepare for applications to $\Ginf$-regularity theory for pseudodifferential equations with generalized symbols, one may, as an intermediate step, investigate the propagation of  $\Ginf$-singularities in case of differential equations with smooth coefficients. Having this special situation in mind, it is natural to consider the following set, defined for any $u\in\G(\Om)$ by
\beq
\label{characteristic_inter}
\Wcl(u) := \bigcap\Char(A),
\eeq
where the intersection is taken over all classical properly supported operators $A\in\Psi^0(\Om)$ such that $Au\in\Ginf(\Om)$. A careful inspection of the proof of Theorem \ref{theorem_wave_front} suggests that $\Wcl(u)$ can be used in place of $\Wsc(u)$. Indeed, it gives an alternative characterization of $\WFg(u)$.
\begin{theorem}
\label{theorem_Wcl}
For all $u\in\G(\Om)$
\beq
\label{Wcl=Wsc=WFg}
\Wcl(u) = \Wsc(u) = \WFg(u).
\eeq  
\end{theorem}
\begin{proof}
The inclusion $\Wsc(u)\subseteq\Wcl(u)$ is obvious. Conversely, if $(x_0,\xi_0)\not\in{\WFg(u)}$, as in the proof of Theorem \ref{theorem_wave_front} one can find a properly supported operator $P\in\Psi^0(\Om)$ such that $(x_0,\xi_0)\notin\Char(P)$ and $Pu\in\Ginf(\Om)$.
\end{proof}

\section{Noncharacteristic $\Ginf$-regularity and\\
propagation of singularities}

As a first application of Theorem \ref{theorem_wave_front} we prove an
extension of the classical result on noncharacteristic regularity for
distributional solutions of arbitrary pseudodifferential equations (with
smooth symbols). A generalization of this result for partial differential
operators with Colombeau coefficients was achieved in \cite{HOP:03}, here we
present a version for pseudodifferential operators with slow scale generalized
symbols. 
\begin{theorem} If $P = p(x,D)$ is a properly supported pseudodifferential
operator 
with slow scale symbol and $u\in\G(\Om)$ then \beq\label{nonchar_theorem}
  \WFg(P u) \subseteq \WFg(u) \subseteq
    \WFg(P u) \cup \compl{\Ellsc(p)}.
\eeq
\end{theorem}
\begin{proof} The first inclusion relation is obvious from Theorem
\ref{theorem_micro_support}. \\
Assume that $(x_0,\xi_0) \not\in \WFg(P u)$ but $p$ is micro-elliptic there.
Thanks to Theorem \ref{theorem_wave_front} we can find $a(x,D) \in
    \Oprop{0}(\Om)$ such that $a(x,D) P u \in \Ginf(\Om)$ and $(x_0,\xi_0) \in
\Ellsc(a)$. By the (slow scale) symbol calculus  and Proposition
\ref{proposition_lower_order}(ii) $a(x,D) p(x,D)$ has a slow scale symbol,
which is micro-elliptic 
at $(x_0,\xi_0)$. Hence another application of Theorem \ref{theorem_wave_front}
yields that $(x_0,\xi_0) \not\in \WFg(u)$.
\end{proof}

\begin{remark} As can be seen from various examples in \cite{HO:03},
relation (\ref{nonchar_theorem}) does not hold in general for regular symbols
$p$ which satisfy estimate (\ref{estimate_below}). In this sense, the overall
slow scale properties of the symbol are crucial in the above statement and are
not just technical convenience. In fact, adapting the reasoning in
\cite[Example 4.6]{HO:03} to the symbol $p_\eps(x,\xi) = 1 + c_\eps x^2$,
$c_\eps \geq 0$, we obtain the following: $p_\eps(x,\xi) \geq 1$ whereas the
unique solution $u$ of $p u = 1$ is $\Ginf$ if and only if $(c_\eps)_\eps$ is a
slow scale net.
\end{remark}

In the remainder of this section we show how Theorem \ref{theorem_Wcl} enables
us to extend a basic result on propagation of singularities
presented in \cite[Sections 23.1]{Hoermander:V3}, where we now allow for  
Colombeau generalized functions as solutions and initial values in first-order
strictly hyperbolic partial differential equations with smooth coefficients.
Hyperbolicity will be assumed with respect to time direction and we will
occasionally employ pseudodifferential operators whose symbols depend
smoothly on the real parameter $t$. This means that we have symbols from the
space $\Cinf(\R,S^m(\Om\times\R^n))$: these are 
of the form $a(t,x,\xi)$ with $a\in\Cinf(\R\times\Om\times\R^n)$ such that  
for each $t$ and $h\in\N$ one has $\frac{d^{h}}{dt^h}a(t,.,.)\in S^m(\Om\times\R^n)$, where all symbol
seminorm estimates are uniform when $t$ varies in a compact subset
of $\R$. Defined on representatives in the obvious way, a properly supported 
operator $a(t,x,D_x)$ maps $\G(\Om)$ into $\G(\Om)$. Moreover, if we assume that $a(t,x,D_x)$ is uniformly properly supported with respect to $t$, that is there exists a proper closed set $L$ such that $\text{supp}\,k_{a(t,x,D_x)}\subseteq L$ for all $t$, if $u\in\G(\Om\times\R)$ then $a(t,x,D_x)(u(t,\cdot))\in\G(\Om\times\R)$.
\begin{proposition}\label{prop_first_order}
Let $P(t,x,D_x)$ be a first-order partial differential operator with 
real principal symbol $P_1$ and coefficients in $\Cinf(\R\times\Om)$, which
are constant outside some compact subset of $\Om$.  
Assume that $u\in\G(\Om\times\R)$ satisfies the homogeneous Cauchy problem
 \begin{align}
 \d_t u + i P(t,x,D_x) u &= 0 
		\label{first_order_equ}\\
 u(.,0) &= g \in\G(\Om) \label{initial_cond}.
 \end{align}
If $\Phi_t$ denotes the Hamilton flow corresponding to $P_1(t,.,.)$ on
$T^*(\Om)$ then we have for all $t\in\R$  
 \beq\label{WFg_flow}
    \WFg(u(.,t)) = \Phi_t\big( \WFg(g) \big).
 \eeq
\end{proposition}
\begin{proof}
We adapt the symbol constructions presented in \cite{Hoermander:V3}, pp.\
388-389. Observe that one can carry out all steps of that classical procedure
in $\Om \subseteq \R^n$ and with all operators uniformly properly supported. To be
more precise, let $(x_0,\xi_0)\in\big(\CO{\Om}\big)\setminus\WFg(g)$ and
choose $q\in S^0(\Om\times\R^n)$ polyhomogeneous, i.e., having homogeneous
terms in the asymptotic expansion, such that $q$ is micro-elliptic at
$(x_0,\xi_0)$, $\mu\supp(q)\cap\WFg(g)= \emptyset$, and $q(x,D)$ is properly
supported. By Corollary \ref{corollary_wave_front} we deduce that
$\WFg(q(x,D)g) \subseteq \WFg(g) \cap \mu\supp(q) = \emptyset$, therefore
$q(x,D)g\in\Ginf(\Om)$.

We can find a symbol $Q(t,x,\xi)$, which is polyhomogeneous of order $0$,
smoothly depending on the parameter $t\in\R$, and with the following 
properties: 
\begin{trivlist}
\item the operators $Q(t,x,D_x)$ are uniformly properly supported for $t\in\R$,
\item $[\d_t+iP(t,x,D_x),Q(t,x,D_x)]=R(t,x,D_x)$ is a $t$-parametrized
	operator\\  of order $-\infty$ and uniformly properly supported,
\item $Q_0(t,x,\xi)=q_0(\Phi_t^{-1}(x,\xi))$ for the principal symbols,
\item and $Q(0,x,\xi)-q(x,\xi)\in S^{-\infty}(\Om\times\R^n)$.
\end{trivlist}
From the properties of $Q$ and \eqref{first_order_equ} we have 
\begin{multline*}
(\d_t+i P(t,x,D_x))Q(t,x,D_x)u \\
	= Q(t,x,D_x)(\d_tu+i P(t,x,D_x)u)+[\d_t+iP,Q]u 
= R(t,x,D_x)u,
\end{multline*}
and $Q(0,x,D_x)u(.,0)=q(x,D)g+R_0(x,D)g$, where $R_0$ is of order $-\infty$. Observe that $R(t,x,D_x)u(.,t) \in
\Ginf(\Om)$ and $R_0(x,D)g\in\Ginf(\Om)$, which implies
$Q(0,x,D_x)u(.,0)\in\Ginf(\Om)$. Setting $v=Qu\in\G(\Om\times\R)$ we obtain $\d_t v + i P(t,x,D_x) v\in\G(\Om\times\R)$ and
\begin{align*} 
 \d_t v(.,t) + i P(t,x,D_x) v(.,t) & \in\Ginf(\Om) \qquad \forall t\in\R,\\ 
	v(.,0) & \in \Ginf(\Om), 
 \end{align*} 
which we interpret as a Cauchy problem with
$\Ginf$-data with respect to the space variable $x$. From an inspection of the
energy estimates discussed in \cite{LO:91} one directly infers that
$v(.,t) \in \Ginf(\Om)$ for all $t\in\R$ (note that the coefficients are
independent of 
$\epsilon$). Therefore at fixed value of $t$ we have that
$Q(t,x,D_x)u(.,t)\in\Ginf(\Om)$ and the noncharacteristic regularity relation
\eqref{nonchar_theorem} yields
 \begin{multline}\label{WF_u_Q}
	\WFg(u(.,t))\subseteq \Char(Q(t,x,D_x))\\
	=\{ (x,\xi):\, Q_0(t,.,.)=q_0(\Phi_t^{-1}(.,.)) 
	\text{ is not micro-elliptic at } (x,\xi) \}. 
 \end{multline}
But $q_0\circ\Phi_t^{-1}$ is micro-elliptic at $(x,\xi)=\Phi_t(x_0,\xi_0)$, which implies $(x,\xi)\not\in\WFg(u(.,t))$ and therefore we have shown 
\[
\WFg(u(.,t)) \subseteq \Phi_t(\WFg(g)).
\]
The opposite inclusion is proved by time reversal.
\end{proof}

\begin{remark} The same result can be proven when $P$ is a pseudodifferential
operator (with parameter $t$ and global symbol estimates), but with $g$ and
$u$ in the $\G_{2,2}$-spaces as introduced in \cite{BO:92} and using a basic
regularity result from \cite{GH:03}. However, this formally requires first to
transfer the notion of micro-ellipticity to that context, which will be left
to future presentations.
\end{remark}

While Proposition \ref{prop_first_order} determines the wave front set of
$u(.,t)$ for fixed $t$, we are aiming for a complete description of $\WFg(u)$
in the cotangent bundle over $\Om\times\R$. The crucial new ingredient needed
in extending the analogous discussion in \cite{Hoermander:V3}, p.\ 390, to
Colombeau generalized functions is a result about microlocal regularity of the
restriction operator $\G(\Om\times\R) \to \G(\Om)$.
\begin{lemma}\label{restriction}
Let $I\subseteq\R$ be an open interval, $t_0\in I$, and $u\in\G(\Om\times I)$
such that 
 \beq\label{inter_cond}
  \WFg(u) \cap \big( \Om\times\{t_0\}\times \{(0,\tau) : \tau\in\R \}\big) 
	= \emptyset. 
 \eeq
Then the wave front set of the restriction of $u$ to $\Om\times\{t_0\}$
satisfies the following relation
 \beq\label{WFg_restriction}
    \WFg(u \mid_{t=t_0}) \subseteq \{ (x,\xi)\in\CO{\Om} : \exists\tau\in\R:
    (x,t_0;\xi,\tau)\in\WFg(u)  \}.
 \eeq
\end{lemma}
\begin{remark}
Note that, unlike with distributions, the restriction $u \mid_{t=t_0}$ is
always well-defined for Colombeau functions on $\Om\times I$, regardless of the
microlocal intersection condition (\ref{inter_cond}) in the above lemma.
However, an obvious adaption of the counter example in \cite[Example
5.1]{HK:01} shows that, in general, the latter cannot be dropped without 
losing the relation (\ref{WFg_restriction}).
\end{remark}
\begin{proof}
We are proving a microlocal statement, so we may assume that $t_0 = 0$ and 
$u$ has compact support near $\{ t = 0 \}$; in particular, we may then pick a
representative $(u_\eps)_\eps$ of $u$ with $\supp(u_\eps)$ contained in a
fixed compact set uniformly for all $\eps\in (0,1]$. 

The aim of the proof is to show the following: if $(x_0,\xi_0)\in\CO{\Om}$ is 
in the complement of the right-hand side of
relation (\ref{WFg_restriction}) and $\phi\in\Cinfc(\Om)$ supported near $x_0$,
then $(\phi(.) u_\eps(.,0))\FT{\ }(\xi)$ is (Colombeau-type) rapidly
decreasing, i.e., with uniform $\eps$-asymptotics (cf.\ \cite[Definition
17]{Hoermann:99}), in some conic neighborhood $\Ga(\xi_0)$ of $\xi_0$. As a
preliminary observation, note that if $\psi\in\Cinfc(I)$ with 
$\psi(0) = 1$ then 
we may write
\begin{align}\label{restr_integral}
    \big(\phi(.) u_\eps(.,0)\big)\FT{\ }(\xi) & =
	\big((\phi\otimes\psi)(.,0) u_\eps(.,0)\big)\FT{\ }(\xi) \\
     & = \int \F_{n+1}\big((\phi\otimes\psi)\, u_\eps \big)(\xi,\tau) 
		\,\dslash\tau, \nonumber
\end{align}
where $\F_{n+1}$ denotes $(n+1)$-dimensional Fourier transform. We will find
rapid decrease estimates of the integrand upon an appropriate splitting of 
the integral depending on the parameter $\xi$.

First, the hypothesis (\ref{inter_cond}) gives that for each $y\in\Om$ we find
an open neighborhood $V(y,0)$ and open cones $\Ga^\pm_y \ni (0,\pm1)$ such
that for all $f\in\Cinfc(V(y,0))$ the function $\FT{(f u_\eps)}(\xi,\tau)$ is
Colombeau-rapidly decreasing in $\Ga_y := \Ga^-_y \cup \Ga^+_y$. By
compactness the set $\supp(u) \cap (\Om\times\{0\})$ is covered by 
finitely many neighborhoods $V(y_1,0),\ldots, V(y_M,0)$. Again by compactness,
we may assume that there is $U_1\subseteq \Om$ open and $\eta_1 > 0$ such that
$\supp(u) \cap (\Om\times\{0\})\subseteq U_1\times (-\eta_1,\eta_1) 
\subseteq \bigcup_{j=1}^M V(y_j,0)$. Furthermore, $\bigcap_{j=1}^M \Ga_{y_j}$
is an open cone around $(0,-1)$ and $(0,1)$ and there is $c_1 > 0$ such
that it still contains the conic neighborhood $\Ga_1 := \{ (\xi,\tau) : 
|\tau| \geq c_1 |\xi| \}$. Using a finite partition of unity subordinated to
$(V(y_j,0))_{1 \leq j \leq M}$, we obtain the following statement: for any
$\phi\in\Cinfc(U_1)$ and 
$\psi\in\Cinfc((-\eta_1,\eta_1))$ there exists $N \geq 0$ with the 
property that $\forall l\in\N\, \exists C_l>0\, \exists \eps_0 > 0$ 
which guarantee the rapid decrease estimate
 \beq\label{rap_dec_1}
    |\F_{n+1} \big((\phi\otimes\psi)\, u_\eps\big)(\xi,\tau)|
        \leq C_l \eps^{-N} \lara{(\xi,\tau)}^{-l}  \qquad 
        (\xi,\tau)\in\Ga_1, 0 < \eps < \eps_0.
 \eeq
This will provide corresponding upper bounds for the integrand in
(\ref{restr_integral}) whenever $|\tau| \geq c_1 |\xi|$, $\xi$ arbitrary.

Second, it follows from the assumption on $(x_0,\xi_0)$ that 
$\big(\{(x_0,0)\}\times \{\xi_0\} \times\R\big) \cap \WFg(u) = \emptyset$. 
Hence for all $\sigma\in\R$ there is an open set $V_{\sigma}\ni(x_0,0)$
and an open conic neighborhood
$\Ga(\xi_0,\sigma)$ such that for all $f\in\Cinfc(V_{\sigma})$ the function 
$\FT{(f u_\eps)}(\xi,\tau)$ is Colombeau-rapidly decreasing in $\Ga_{\sigma}$.
By a conic compactness argument (via projections to the
unit sphere), we deduce that finitely many cones $\Ga_{\sigma_j}$
($j=1,\ldots,M$) suffice to cover the two-dimensional sector 
$\{(\la \xi_0,\tau) : \la > 0, |\tau| \leq c_1 \la |\xi_0| \}$. 
Let $\pi_n\colon\R^{n+1} \to \R^n$ be the projection $\pi_n(\xi,\tau) = \xi$
and define 
\[
    \Ga(\xi_0) := \bigcap_{j=1,\ldots,M} \pi_n(\Ga_{\sigma_j}),
\]
which is an open conic neighborhood of $\xi_0$ in $\R^n\zs$. Furthermore, 
$\bigcap_{j=1}^M V_{\sigma_j} \ni (x_0,0)$ is open and contains still some
neighborhood of product form, say $U_0\times (-\eta_0,\eta_0)$. 
Therefore, for any $\phi\in\Cinfc(U_0)$ and
$\psi\in\Cinfc((-\eta_0,\eta_0))$ there exists $N \geq 0$ with the following
property: 
$\forall l\in\N\, \exists C_l>0\, \exists \eps_0 > 0$ we have an 
estimate
 \beq\label{rap_dec_0}
    |\F_{n+1}((\phi\otimes\psi)\, u_\eps)(\xi,\tau)|
        \leq C_l \eps^{-N} \lara{(\xi,\tau)}^{-l} \qquad 
	\xi\in\Ga(\xi_0), |\tau| \leq c_1 |\xi|, 
 \eeq
when $0 < \eps < \eps_0$.
So if $\xi\in\Ga(\xi_0)$, we also obtain 
 corresponding upper bounds in (\ref{restr_integral}) 
over the remaining integration domain $|\tau| \leq c_1 |\xi|$.

To summarize, when $\xi\in \Ga(\xi_0)$ we may combine 
(\ref{rap_dec_1}-\ref{rap_dec_0}) by taking $U := U_0 \cap U_1$,
$\phi\in\Cinfc(U)$, and $\eta := \min(\eta_0,\eta_1)$, $\psi\in\Cinfc((-\eta,\eta))$,
$\psi(0)=1$. Upon applying this to (\ref{restr_integral}) we arrive at
\[
    |(\phi(.) u_\eps(.,0))\FT{\ }(\xi)| \leq  \eps^{-N} C_l
        \int \!\!\! \frac{\dslash\tau}{(1 + |\xi|^2 + \tau^2)^{l/2}}
        = \eps^{-N} C_l \lara{\xi}^{1-l}\!\! 
		\int\limits_\R \!\lara{r}^{-l} \,\dslash r,
\]
for some $N$ independent of $l \geq 2$ and $\eps$ sufficiently small.
\end{proof}

\begin{theorem}\label{thm_first_order}
Let $u$ be the (unique) solution of the homogeneous Cauchy problem
(\ref{first_order_equ}-\ref{initial_cond}) and denote by 
$\ga(x_0,\xi_0)$ the maximal bicharacteristic curve passing through
the characteristic point 
$(x_0,0;\xi_0,-P_1(0,x_0,\xi_0))\in\CO{\Om\times\R}$. Then
the generalized wave front of $u$ is given by
 \beq
    \WFg(u) = \bigcup_{(x_0,\xi_0)\in\WFg(g)} \ga(x_0,\xi_0).
 \eeq
\end{theorem}
\begin{proof} We refer to the discussion in \cite{Hoermander:V3}, p.\ 390,
which we can adapt to our case with only little changes required. 
First, $\d_t + i P(t,x,D_x)$ is already a partial differential operator and
hence we 
obtain $\WF_g(u) \subseteq \Char(\d_t + i P_1) \not\ni (x,t;0,\tau)$ due to
(\ref{nonchar_theorem}); therefore, combined with Lemma \ref{restriction} we
immediately obtain the inclusion 
 \[
	\WFg(u(.,t)) \subseteq \{ (x,\xi)\in\CO{\Om} : 
	(x,t;\xi,-P_1(t,x,\xi)) \in \WFg(u) \}.
 \]

On the other hand, based on the inclusion relation (\ref{WF_u_Q}) we can carry out the following construction: let $(x_1,t_1,\xi_1,\tau_1)\in\CO{\Omega\times\R}$ such that $\xi_1\neq 0$ and $(x_1,\xi_1)\not\in\WFg(u(.,t_1))$ and $Q_0$ as in the proof of Proposition \ref{prop_first_order}, which is micro-elliptic at $(x_1,\xi_1)$. We claim that $(x_1,t_1,\xi_1,\tau_1)\not\in\WFg(u)$. 
Indeed, one may use cut-off functions of product form $\phi(x) \psi(t)$ with
small supports near $x_1$, $t_1$ and write 
$\big((\phi\otimes\psi) u_\eps\big)\FT{\ }(\xi,\tau) = \int e^{-it\tau}\psi(t)(\phi u_\eps(.,t))\FT{\ }(\xi)\, dt$. We have Colombeau rapid decrease estimates in $(\xi,\tau)$ for those 
$\xi$-directions where $Q(t,x,\xi)$ stays micro-elliptic for all $(x,t) \in
	\supp(\phi\otimes\psi)$.  
Together with the observation at the
	beginning of the proof and \eqref{WFg_flow} this implies the equality 
 \begin{multline*}
 \Phi_t(\WFg(g)) = \WFg(u(.,t)) \\= \{ (x,\xi)\in\CO{\Om} : 
		(x,t;\xi,-P_1(t,x,\xi)) \in \WFg(u) \},
 \end{multline*}
which  yields the asserted statement. 
\end{proof}

\section*{Appendix: Pseudodifferential calculus with\\
\hphantom{Appendix:}
general scales}

\setcounter{section}{1}
\renewcommand{\thesection}{\Alph{section}}
\renewcommand{\theequation}{A.\arabic{equation}}
\setcounter{equation}{0}
 
We provide some background of the required pseudodifferential tools for a sufficiently large class of generalized symbols.
\begin{definition}\label{def_general}
Let $\mathfrak{A}$ be the set of all nets $(\omega_\eps)_\eps\in\R^{(0,1]}$ such that $c_0\le\omega_\eps\le c_1\epsilon^{-p}$ for some $c_0,c_1,p > 0$ and for all $\eps$. Let $\mathfrak{B}$ be any subset of $\mathfrak{A}$ closed with respect to pointwise product and maximum. For $m\in\R$ and $\Omega$ an open subset of $\R^N$, we define the spaces of $\mathfrak{B}$-nets of symbols
\[
\begin{split}
 \Sbeu^m(\Om\times\R^n) := \{
(a_\eps)_\eps\in \mathcal{S}^m[\Om\times\R^n]:\, &\forall K \Subset \Omega\ \exists (\omega_\eps)_\eps \in \mathfrak{B}\\
& \forall \alpha, \beta \in \N^n\, \exists c>0\,  \forall \eps:\, |a_\eps|^{(m)}_{K,\alpha,\beta} \le c\, \omega_\eps\},
\end{split}
\]
\[
\begin{split}
 \Sbeu^{-\infty}(\Om\times\R^n) := \{
(a_\eps)_\eps\in \mathcal{S}^m[\Om\times\R^n]&:\, \forall K \Subset \Omega\ \exists (\omega_\eps)_\eps \in \mathfrak{B}\ \forall m \in \R\\
& \forall \alpha, \beta \in \N^n\, \exists c>0\,  \forall \eps:\, |a_\eps|^{(m)}_{K,\alpha,\beta} \le c\, \omega_\eps\}
\end{split}
\]
and the factor spaces
\[
 \Sybeu^m(\Om\times\R^n) := \Sbeu^m(\Om\times\R^n) / \Nu^m(\Om\times\R^n),
\]
\[
 \Sybeu^{m,-\infty}(\Om\times\R^n) := \Sbeu^m(\Om\times\R^n) / \Nuinf(\Om\times\R^n)
\]
and 
\[
\Sybeu^{-\infty}(\Om\times\R^n) := \Sbeu^{-\infty}(\Om\times\R^n) / \Nuinf(\Om\times\R^n).
\]
If\/ $\Omega$ is an open subset of $\R^n$ then the elements of $\Sbeu^m(\Om\times\R^n)$, $\Sybeu^m(\Om\times\R^n)$, and $\Sybeu^{m,-\infty}(\Om\times\R^n)$ are called \emph{$\mathfrak{B}$-nets of symbols of order $m$}, \emph{$\mathfrak{B}$-generalized symbols of order $m$}, and \emph{$\mathfrak{B}$-generalized symbols of refined order $m$} respectively. The sets $\Sbeu^{-\infty}(\Om\times\R^n)$ and $\Sybeu^{-\infty}(\Om\times\R^n)$ constitute the \emph{$\mathfrak{B}$-nets of smoothing symbols} and  the \emph{$\mathfrak{B}$-generalized smoothing symbols} respectively. Finally, if instead of $\Om$ we have $\Om\times\Om$ then we use the notion of amplitude rather than symbol.  
\end{definition}
As already mentioned in \cite{GGO:03}, if we require the above symbol estimates
to hold only for small values of $\eps$
this would result in larger spaces of nets of symbols as well as somewhat
larger quotient spaces. Even though it is possible then to define
pseudodifferential operators with $\mathfrak{B}$-generalized symbols equally
well
we prefer to consider here the spaces $\Sybeu^m(\Om\times\R^n)$ defined above
since a complete pseudodifferential calculus can be developed.\\

\begin{remark} 
\label{rem_special}
Note that as special cases of $\mathfrak{B}$-generalized
symbols we obtain for $\mathfrak{B}=\Pi_\ssc$ the slow scale generalized
symbols, and for $\mathfrak{B}=\{(\eps^{-N})_\eps : N\in\N\}$ the regular
generalized symbols $\Syru^m(\Om\times\R^n)$ introduced in earlier work (cf.\ \cite{GGO:03}), to which
we refer for further details and notations concerning regular symbols and
amplitudes.
\end{remark}

In the sequel we say that $(a_\eps)_\eps\in\Sbeu^m(\Om\times\R^n)$ is of growth
type $(\omega_\eps)_\eps\in\mathfrak{B}$ on $K \Subset\Omega$ if
$(\omega_\eps)_\eps$ estimates each seminorm $|a_\eps|^{(m)}_{K,\alpha,\beta}$.
We list the basic steps in establishing a calculus with asymptotic
expansions.\begin{definition}
\label{def_asym_expan_represen}
Let $\{m_j\}_j$ be a decreasing sequence of real numbers tending to $-\infty$ and let $\{(a_{j,\eps})_\eps\}_j$ be a sequence of $\mathfrak{B}$-nets of symbols $(a_{j,\eps})_\eps\in\Sbeu^{m_j}(\Om\times\R^n)$ such that
\beq\label{growth_type}
 \forall K \Subset \Om\ \exists (\omega_\eps)_\eps \in \mathfrak{B}\ \forall j\in\N:\ (a_{j,\eps})_\eps\, \text{is of growth type $(\omega_\eps)_\eps$ on $K$}.
\eeq
We say that $\sum_j(a_\eps)_\eps$ is the asymptotic expansion of $(a_\eps)_\eps\in\E[\Om\times\R^n]$, $(a_\eps)_\eps\sim\sum_j(a_{j,\eps})_\eps$ for short, iff for all $K\Subset\Omega$ there exists $(\omega_\eps)_\eps\in\mathfrak{B}$ such that for all $r\ge 1$ the difference $(a_\eps-\sum_{j=0}^{r-1}a_{j,\eps})_\eps$ belongs to $\Sbeu^{m_r}(\Om\times\R^n)$ and is of growth type $(\omega_\eps)_\eps$ on $K$.
\end{definition}
\begin{theorem}
\label{theorem_asym_expan_represen}
For any sequence of $\mathfrak{B}$-nets of symbols $(a_{j,\eps})_\eps\in\Sbeu^{m_j}(\Om\times\R^n)$ as in Definition \ref{def_asym_expan_represen} there exists $(a_\eps)_\eps\in\Sbeu^{m_0}(\Om\times\R^n)$ such that $(a_\eps)_\eps\sim\sum_j(a_{j,\eps})_\eps$. Moreover if $(a'_\eps)_\eps\sim\sum_j(a_{j,\eps})_\eps$ then $(a_\eps-a'_\eps)\in\Sbeu^{-\infty}(\Om\times\R^n)$.
\end{theorem}
This result is easily obtained from the proof of Theorem 5.3 in \cite{GGO:03}, noting we do not have powers of $\eps$ depending on $x$-derivatives and replacing $\eps^{-N}$ with $(\omega_\eps)_\eps\in\mathfrak{B}$. The following proposition concerning negligible nets of symbols is a consequence of Theorem 5.4 in \cite{GGO:03}. 
\begin{proposition}
\label{proposition_asym_expan_negligible}
Let $\{m_j\}_j$ be a decreasing sequence of real numbers tending to $-\infty$ and let $(a_{j,\eps})_\eps\in\Nu^{m_j}(\Om\times\R^n)$ for all $j$. Then there exists $(a_\eps)_\eps\in\Nu^{m_0}(\Om\times\R^n)$ such that for all $r\ge 1$ we have  $(a_\eps-\sum_{j=0}^{r-1}a_{j,\eps})_\eps\in\Nu^{m_r}(\Om\times\R^n)$.
\end{proposition}
\begin{definition}
\label{def_asym_expan_symbol}
Let $\{m_j\}_j$ be a decreasing sequence of real numbers tending to $-\infty$. Let $\{a_j\}_j$ be a sequence of $\mathfrak{B}$-generalized symbols $a_j\in\Sybeu^{m_j}(\Om\times\R^n)$ such that there exists a choice of representatives $(a_{j,\eps})_\eps$ of $a_j$ satisfying \eqref{growth_type}. We say that $\sum_j a_j$ is the asymptotic expansion of $a\in\Sybeu^{m_0}(\Om\times\R^n)$, $a\sim\sum_j a_j$ for short, iff there exists a representative $(a_\eps)_\eps$ of $a$ and for all $j$ representatives $(a_{j,\eps})_\eps$ of $a_j$, such that $(a_\eps)_\eps\sim\sum_j(a_{j,\eps})_\eps$.
\end{definition}
Proposition \ref{proposition_asym_expan_negligible} allows us to claim that $a\sim\sum_j a_j$ iff for any choice of representatives $(a_{j,\eps})_\eps$ of $a_j$ there exists a representative $(a_\eps)_\eps$ of $a$ such that $(a_\eps)_\eps\sim\sum(a_{j,\eps})_\eps$. This observation combined with Theorem \ref{theorem_asym_expan_represen} is crucial in the proof of Theorem \ref{theorem_asym_expan_symbol}.
\begin{theorem}
\label{theorem_asym_expan_symbol}
For any sequence of $\mathfrak{B}$-generalized symbols $a_j\in\Sybeu^{m_j}(\Om\times\R^n)$ as in Definition \ref{def_asym_expan_symbol} there exists $a\in\Sybeu^{m_0}(\Om\times\R^n)$ such that $a\sim\sum_j a_j$. Moreover, if $b\in\Sybeu^{m_0}(\Om\times\R^n)$ has asymptotic expansion $\sum_j a_j$ then there exists a representative $(a_\eps)_\eps$ of $a$ and a representative $(b_\eps)_\eps$ of $b$ such that $(a_\eps-b_\eps)_\eps\in\Sbeu^{-\infty}(\Om\times\R^n)$.
\end{theorem}
\begin{remark}
\label{remark_negligible_infty}
Definition \ref{def_asym_expan_symbol} can be stated for $\mathfrak{B}$-generalized symbols of refined order. More precisely, if $\{m_j\}_j$ is a sequence as above, $a\in\Sybeu^{m_0,-\infty}(\Om\times\R^n)$, $a_j\in\Sybeu^{m_j,-\infty}(\Om\times\R^n)$ for all $j$ and $(a_\eps)_\eps$ and $(a_{j,\eps})_\eps$ denote representatives of $a$ and $a_j$ respectively, the following assertions are equivalent:
\begin{itemize}
\item[(1)] $a \sim \sum_j a_j$,
\item[(2)] $\exists (a_\eps)_\eps$\, $\exists \{(a_{j,\eps})_\eps\}_j\, :$ $(a_\eps)_\eps \sim \sum_j(a_{j,\eps})_\eps$,
\item[(3)] $\forall \{(a_{j,\eps})_\eps\}_j$\, $\exists (a_\eps)_\eps\, :$ $(a_\eps)_\eps \sim \sum_j(a_{j,\eps})_\eps$,
\item[(4)] $\forall \{(a_{j,\eps})_\eps\}_j$\, $\forall (a_\eps)_\eps\, :$ $(a_\eps)_\eps \sim \sum_j(a_{j,\eps})_\eps$.
\end{itemize}
\end{remark}
We briefly recall the main definitions and results concerning pseudodifferential operators with $\mathfrak{B}$-generalized symbols. As already observed after Theorem \ref{theorem_asym_expan_represen}, the proofs are obtained from the corresponding ones in \cite{GGO:03}, with the slight difference of having a net $(\omega_\eps)_\eps$ in $\mathfrak{B}$ instead of a power $\eps^{-N}$.
\begin{definition}
\label{def_pseudo}
Let $a\in\Sybeu^{m}(\Om\times\R^n)$. The pseudodifferential operator with $\mathfrak{B}$-generalized symbol $a$, is the map $a(x,D):\Gc(\Om)\to\G(\Om)$ given by the formula
\[
 a(x,D)u := \int_{\R^n}e^{ix\xi}a(x,\xi)\widehat{u}(\xi)\, \dslash\xi := \biggl[ \biggl( \int_{\R^n} e^{ix\xi} a_\eps(x,\xi) \widehat{u}_\eps(\xi) \dslash\xi \biggr)_\eps \biggr] .
\]
\end{definition}
Proposition 4.7 in \cite{GGO:03} guarantees the well-definedness of $a(x,D)$ as well as the additional mapping property $a(x,D):{\Gcinf}(\Om)\to\Ginf(\Om)$.

In the following, we make
occasional use of some basic properties of the  
space $L(\Gc(\Om),\wt\C)$ of $\wt\C$-linear maps from $\Gc(\Om)$ to $\wt\C$.
In particular, we recall that $\G(\Om)$ is linearly embedded into
$L(\Gc(\Om),\wt\C)$ via generalized integration and $L(\Gc(\Om),\wt\C)$ is a sheaf with respect to $\Om$. This and further results 
are discussed in detail in \cite[Section 2]{GGO:03}.

\begin{definition}
\label{def_kernel} 
Let $a\in\Sybeu^m(\Om\times\R^n)$. The kernel of $a(x,D)$ is the 
 $\wt{\C}$-linear map $k \colon \Gc(\Om\times\Om) \to \wt{\C}$
defined by
\beq
\label{kernel}
k(u):=\int_\Om a(x,D)(u(x,\cdot))\, dx .
\eeq
\end{definition}
To see that formula \eqref{kernel} makes sense for $k$ as an element of 
$L(\Gc(\Om\times\Om),\wt\C)$, we may reason as in 
\cite[Proposition 3.10 and Remark 3.11]{GGO:03} that we have 
$a(x,D)(u(x,\cdot)) \in \Gc(\Om)$. Moreover, for all $u,v\in\Gc(\Om)$
\[
 k(v\otimes u) = \int_\Om a(x,D)u\, v(x)\, dx = \int_\Om u(x){\ }^ta(x,D)v\, dx ,
\]
where $v\otimes u := [(v_\eps(x)u_\eps(y))_\eps] \in \Gc(\Om\times\Om)$; as a
consequence, since $\G(\Om)$ is embedded into $L(\Gc(\Om),\wt\C)$, pseudodifferential operators having the same kernel are identical.\\
We say that a pseudodifferential operator with $\mathfrak{B}$-generalized symbol is \emph{properly supported} if the support of its kernel is a proper set of $\Om\times\Om$. As shown in Proposition 4.17 of \cite{GGO:03}, we have that any properly supported pseudodifferential operator $a(x,D)$ maps $\Gc(\Om)$ into $\Gc(\Om)$, $\Gcinf(\Om)$ into $\Gcinf(\Om)$ and can be extended uniquely to a linear map from $\G(\Om)$ into $\G(\Om)$ such that for all $u\in\G(\Om)$ and $v\in\Gc(\Om)$
\[
 \int_\Om a(x,D)u\, v(x)\, dx = \int_\Om u(x){\ }^ta(x,D)v\, dx .
\]
This extension maps $\Ginf(\Om)$ into $\Ginf(\Om)$. By the same reasoning as in Proposition 4.11 in \cite{GGO:03}, we prove that each pseudodifferential operator $a(x,D)$ with $\mathfrak{B}$-generalized symbol has the pseudolocality property, i.e., $$\singsupp_g(a(x,D)u) \subseteq \singsupp_g(u)$$ for all $u\in\Gc(\Om)$, and that this result is valid for all $u$ in $\G(\Om)$ if $a(x,D)$ is properly supported.\\
Pseudodifferential operators can be defined also by $\mathfrak{B}$-generalized amplitudes. This means for $b\in\Sybeu^m(\Om\times\Om\times\R^n)$ to define the action of the corresponding operator on $u\in\Gc(\Om)$, via the oscillatory integral
\[
 Bu(x) := \int_{\Om\times\R^n} e^{i(x-y)\xi}b(x,y,\xi)u(y)\, dy\, \dslash\xi ,
\]
which gives a Colombeau function in $\G(\Om)$ (cf. Section 3 in \cite{GGO:03}). It is clear that the same constructions concerning kernel and properly supported pseudodifferential operators are still valid. For the sake of completeness we recall that any \emph{integral operator $R$ with regular kernel}, i.e. any operator of the form
\[
 \qquad Ru(x) = \int_\Om k(x,y)u(y)\, dy,\qquad\qquad\quad u\in\Gc(\Om),
\]
where $k\in\Ginf(\Om\times\Om)$ can be written as a pseudodifferential operator with regular amplitude in $\Syr^{-\infty}(\Om\times\Om\times\R^n)$ vice versa if $B$ is a pseudodifferential operator with amplitude in $\Sybeu^{-\infty}(\Om\times\Om\times\R^n)$, then its kernel is a regular generalized function. Finally an operator with regular kernel is regularizing, i.e. it maps $\Gc(\Om)$ into $\Ginf(\Om)$.\\
Consider $a\in\Sybeu^m(\Om\times\R^n)$, $k$ the kernel of $a(x,D)$ and
let $\chi\in\Cinf(\Om\times\Om)$ be a proper function identically $1$ in a
neighborhood of $\supp\, k$. We may write $a(x,D)=a_0(x,D)+a_1(x,D)$, where 
\beq
\label{equation1}
 a_0(x,D)u:=\int_{\Om\times\R^n} e^{i(x-y)\xi}a(x,\xi)\chi(x,y)u(y)\, dy\, \dslash\xi
\eeq
is a properly supported pseudodifferential operator with generalized amplitude $a(x,\xi)\chi(x,y)\in\Sybeu^m(\Om\times\Om\times\R^n)$ and 
\beq
\label{equation2}
  a_1(x,D)u:=\int_{\Om\times\R^n} e^{i(x-y)\xi}a(x,\xi)(1-\chi(x,y))u(y)\, dy\, \dslash\xi
\eeq
is an operator with regular kernel in $\Ginf(\Om\times\Om)$. The following theorem shows that every properly supported pseudodifferential operator defined via an amplitude can be written in the form of Definition \ref{def_pseudo}. This is the main tool in the proof of Theorem \ref{theorem_calculus}.
\begin{theorem}
\label{theorem_pre_calculus}
For any properly supported pseudodifferential operator $A$ with amplitude $a\in\Sybeu^m(\Om\times\Om\times\R^n)$ there exists $\sigma\in\Sybeu^m(\Om\times\R^n)$ such that $A\equiv \sigma(x,D)$ on $\Gc(\Om)$ and $\sigma \sim \sum_\gamma \frac{1}{\gamma !}\partial^\gamma_\xi D^\gamma_y a(x,y,\xi)_{\vert_{x=y}}$.
\end{theorem}
\begin{theorem}
\label{theorem_calculus}
Let $a\in\Sybeu^m(\Om\times\R^n)$ and $b\in\Sybeu^{m'}(\Om\times\R^n)$ be $\mathfrak{B}$-generalized symbols. If the corresponding pseudodifferential operators are properly supported then
\begin{itemize}
\item[(i)] there exists $a'\in\Sybeu^m(\Om\times\R^n)$ such that ${\ }^ta(x,D)\equiv a'(x,D)$ on $\Gc(\Om)$ and $a'\sim\sum_\gamma \frac{(-1)^{|\gamma|}}{\gamma !}\partial^\gamma_\xi D^\gamma_x a(x,-\xi)$;
\item[(ii)] there exists $a^\ast\in\Sybeu^m(\Om\times\R^n)$ such that $a(x,D)^\ast\equiv a^\ast(x,D)$ on $\Gc(\Om)$ and $a^\ast\sim\sum_\gamma\frac{1}{\gamma !}\partial^\gamma_\xi D^\gamma_x \overline{a}$;
\item[(iii)] there exists $a\sharp b\in\Sybeu^{m+m'}(\Om\times\R^n)$ such that $a(x,D)\circ b(x,D)\equiv a\sharp b(x,D)$ on $\Gc(\Om)$ and $a\sharp b\sim\sum_\gamma \frac{1}{\gamma!}\partial^\gamma_\xi a\, D^\gamma_x b$.
\end{itemize}
\end{theorem}
\begin{proof}
We briefly sketch the proof of assertion $(ii)$. For the details concerning the transposed operator and the product we refer to \cite{GGO:03}. By definition of the formal adjoint, $\langle a(x,D)v,\overline{u}\rangle = \langle v,\overline{a(x,D)^\ast u}\rangle$ for all $u,v\in\Gc(\Om)$. This means \[
 \langle v,\overline{a(x,D)^\ast u} \rangle = \int_\Om \overline{\int_{\Omega\times\R^n}e^{i(x-y)\xi}\overline{a}(y,\xi)u(y)\, dy\, \dslash\xi}\, v(x)\, dx,
\]
which, from the embedding of $\G(\Om)$ into $L(\Gc(\Om),\wt{\C})$, leads us to $\displaystyle a(x,D)^\ast u = \int_{\Om\times\R^n} e^{i(x-y)\xi}\overline{a}(y,\xi)u(y)\, dy\, \dslash\xi.$ Now, $a(x,D)^\ast$ is a properly supported pseudodifferential operator with amplitude $\overline{a}(y,\xi)\in\Sybeu^m(\Om\times\Om\times\R^n)$ and an application of Theorem \ref{theorem_pre_calculus} completes the proof.
\end{proof}
Along the lines of Proposition 5.17 in \cite{GGO:03} we easily prove that the composition of a properly supported pseudodifferential operator with $\mathfrak{B}$-generalized symbol and an operator with regular kernel is an operator with regular kernel. Therefore, combining \eqref{equation1} and \eqref{equation2} with Theorem \ref{theorem_calculus}, we have that for arbitrary pseudodifferential operators with $\mathfrak{B}$-generalized symbol the equalities $(i)$ and $(ii)$ on $\Gc(\Om)$ are valid modulo some operator with regular kernel. Furthermore the composition $a(x,D)\circ b(x,D)$, where at least one of the operators is properly supported, is a pseudodifferential operator $a\sharp b(x,D)$ modulo some operator with regular kernel.
\begin{remark}\leavevmode
\label{remark_rho_delta}
\begin{trivlist}
\item[(i)] It is clear from the structure of $\Sybeu^{m,-\infty}(\Om\times\R^n)$ that all the definitions and results of this appendix can be stated for symbols of refined order.
\item[(ii)] Let now ${\mathcal{S}}^m_{\rho,\delta}[\Om\times\R^n]$ be the set of all nets $(a_\eps)_\eps\in{S}^m_{\rho,\delta}(\Om\times\R^n)^{(0,1]}$ with
\[
|a_\eps|^{(m)}_{K,\alpha,\beta,\rho,\delta} := \sup_{x\in K,\xi\in\R^n}|\partial^\alpha_\xi\partial^\beta_x a_\eps(x,\xi)|\lara{\xi}^{-m+\rho|\alpha|-\delta|\beta|}
\]
seminorm in ${S}^m_{\rho,\delta}(\Om\times\R^n)$. We define $\underline{\mathcal{S}}^m_{\mathfrak{B},\rho,\delta}(\Om\times\R^n)$ and $\underline{\mathcal{N}}^m_{\rho,\delta}(\Om\times\R^n)$ as the subspaces of ${\mathcal{S}}^m_{\rho,\delta}[\Om\times\R^n]$ obtained by requiring the same estimate of $|a_\eps|^{(m)}_{K,\alpha,\beta,\rho,\delta}$ as in $\Sbeu^m(\Om\times\R^n)$ and $\Nu^m(\Om\times\R^n)$ respectively. In conclusion, under the assumption $0 \le \delta < \rho \le 1$ it is possible to develop a pseudodifferential calculus for generalized symbols in ${{\wt{\underline{\mathcal{S}}}}^m_{\mathfrak{B},\rho,\delta}}(\Om\times\R^n) := \underline{\mathcal{S}}^m_{\mathfrak{B},\rho,\delta}(\Om\times\R^n) / \underline{\mathcal{N}}^m_{\rho,\delta}(\Om\times\R^n)$ and ${{\wt{\underline{\mathcal{S}}}}^{m,-\infty}_{\mathfrak{B},\rho,\delta}}(\Om\times\R^n) := \underline{\mathcal{S}}^m_{\mathfrak{B},\rho,\delta}(\Om\times\R^n) / \underline{\mathcal{N}}^{-\infty}(\Om\times\R^n)$, as in the classical theory.
\end{trivlist}
\end{remark}

\paragraph{Acknowledgement:} The authors are grateful to {\sc DIANA}'s
members (cf.\\ {\tt http://www.mat.univie.ac.at/\~{ }diana}) for several
inspiring discussions on the subject. Our special thanks go to Roland
Steinbauer for having posed the striking question leading us to the additional
characterization Theorem \ref{theorem_Wcl}.

\bibliographystyle{abbrv}
\bibliography{claudia}

\newcommand{\SortNoop}[1]{}
\begin{thebibliography}{10}

\bibitem{BO:92}
H.~Biagioni and M.~Oberguggenberger.
\newblock Generalized solutions to the korteweg-de vries and the regularized
  long-wave equations.
\newblock {\em SIAM J. Math Anal.}, 23(4):923--940, 1992.

\bibitem{Colombeau:84}
J.~F. Colombeau.
\newblock {\em New Generalized Functions and Multiplication of Distributions}.
\newblock North Holland, Amsterdam, 1984.

\bibitem{Colombeau:85}
J.~F. Colombeau.
\newblock {\em Elementary Introduction to New Generalized Functions}.
\newblock North-Holland Mathematics Studies 113. Elsevier Science Publishers,
  1985.

\bibitem{DPS:98}
N.~Dapi{\'{c}}, S.~Pilipovi{\'{c}}, and D.~Scarpal{\'{e}}zos.
\newblock Microlocal analysis of {C}olombeau's generalized functions:
  propagation of singularities.
\newblock {\em Jour. d'Analyse Math.}, 75:51--66, 1998.

\bibitem{Folland:95}
G.~B. Folland.
\newblock {\em Introduction to partial differential equations}.
\newblock Princeton University Press, Princeton, New Jersey, second edition,
  1995.

\bibitem{Garetto:04}
C.~Garetto.
\newblock Pseudo-differential operators in algebras of generalized functions
  and global hypoellipticity.
\newblock {\em to appear in Acta Applicandae Mathematicae}, 2004.

\bibitem{GGO:03}
C.~Garetto, T.~Gramchev, and M.~Oberguggenberger.
\newblock Pseudo-differential operators with generalized symbols and regularity
  theory.
\newblock {\em Preprint, University of Innsbruck}, 2003.

\bibitem{GKOS:01}
M.~Grosser, M.~Kunzinger, M.~Oberguggenberger, and R.~Steinbauer.
\newblock {\em Geometric theory of generalized functions}, volume 537 of {\em
  Mathematics and its Applications}.
\newblock Kluwer, Dordrecht, 2001.

\bibitem{Hoermander:71}
L.~H{\"o}rmander.
\newblock {F}ourier integral operators {I}.
\newblock {\em Acta Math.}, 127:79--183, 1971.

\bibitem{Hoermander:V3}
L.~H{\"o}rmander.
\newblock {\em The Analysis of Linear Partial Differential Operators}, volume
  III.
\newblock Springer-Verlag, 1985.
\newblock Second printing 1994.

\bibitem{Hoermann:99}
G.~H{\"o}rmann.
\newblock Integration and microlocal analysis in {C}olombeau algebras.
\newblock {\em J. Math. Anal. Appl.}, 239:332--348, 1999.

\bibitem{GH:03}
G.~H\"{o}rmann.
\newblock First-order hyperbolic pseudodifferential equations with generalized
  symbols.
\newblock {\em to appear in J. Math. Anal. Appl.}, 2004.
\newblock arXiv:math.AP/0307406.

\bibitem{HdH:01}
G.~H{\"o}rmann and M.~V. de~Hoop.
\newblock Microlocal analysis and global solutions of some hyperbolic equations
  with discontinuous coefficients.
\newblock {\em Acta Appl. Math.}, 67:173--224, 2001.

\bibitem{HK:01}
G.~H{\"o}rmann and M.~Kunzinger.
\newblock Microlocal analysis of basic operations in {C}olombeau algebras.
\newblock {\em J. Math. Anal. Appl.}, 261:254--270, 2001.

\bibitem{HO:03}
G.~H{\"o}rmann and M.~Oberguggenberger.
\newblock Elliptic regularity and solvability for partial differential
  equations with {C}olombeau coefficients.
\newblock {\em Electron. J. Diff. Eqns.}, 2004(14):1--30, 2004.

\bibitem{HOP:03}
G.~H\"{o}rmann, M.~Oberguggenberger, and S.~Pilipovic.
\newblock Microlocal hypoellipticity of linear partial differential operators
  with generalized functions as coefficients.
\newblock {\em arXiv:math.AP/0303248}, 2003.

\bibitem{LO:91}
F.~Lafon and M.~Oberguggenberger.
\newblock Generalized solutions to symmetric hyperbolic systems with
  discontinuous coefficients: the multidimensional case.
\newblock {\em J. Math. Anal. Appl.}, 160:93--106, 1991.

\bibitem{NP:98}
M.~Nedeljkov and S.~Pilipovi{\'{c}}.
\newblock Hypoelliptic differential operators with generalized constant
  coefficients.
\newblock {\em Proc. Edinb. Math. Soc.}, 41:47--60, 1998.

\bibitem{NPS:98}
M.~Nedeljkov, S.~Pilipovi{\'{c}}, and D.~Scarpal{\'{e}}zos.
\newblock {\em The Linear Theory of Colombeau Generalized Functions}.
\newblock Pitman Research Notes in Mathematics 385. Longman Scientific {\&}
  Technical, 1998.

\bibitem{O:92}
M.~Oberguggenberger.
\newblock {\em Multiplication of Distributions and Applications to Partial
  Differential Equations}.
\newblock Pitman Research Notes in Mathematics 259. Longman Scientific {\&}
  Technical, 1992.

\end{thebibliography}

\end{document}